\documentclass[microtype]{gtpart}   
\usepackage{pinlabel}
\usepackage[all]{xy}

\author{Brett Parker} 
\givenname{Brett}
\surname{Parker}
\email{brettdparker@gmail.com}
\address{Mathematical Sciences Institute\\
Australian National University\\\newline
Canberra, 2601\\
Australia}


\arxivreference{1102.0158}  
\arxivpassword{kb8eh}   

\usepackage{epsf}  
\usepackage{amsfonts, amssymb, amsthm, amsmath}  
\usepackage{hyperref}
\usepackage{tikz}
\usetikzlibrary{matrix,arrows,decorations.pathmorphing}
\usepackage{tikz-cd}
\usepackage{graphicx}
\usepackage{psfrag}
\newtheorem{thm}{Theorem}  

\newtheorem{cor}[thm]{Corollary}  
\newtheorem{lemma}[thm]{Lemma}  
\newtheorem{remark}[thm]{Remark}  
\newtheorem{defn}[thm]{Definition}  
\newtheorem{prop}[thm]{Proposition}  
\newtheorem{claim}[thm]{Claim}  
\newtheorem{example}[thm]{Example}  
\numberwithin{thm}{section}  
\def\pf{\noindent\emph{Proof: }}  
\def\stop{\hfill$\square$}

\providecommand{\totl}[1]{\ensuremath{\lceil #1\rceil }}
\providecommand{\totb}[1]{\ensuremath{\underline{ #1}}}
\newcommand{\ro}{{}^{r}\Omega}
\newcommand{\rh}{{}^{r}H}

\newcommand{\ex}{\bold}
\providecommand {\e}[1]{\mathfrak t^{#1}}

\providecommand{\Cs}[2]{\ensuremath {C^{#1,\underline{#2}}}}

\providecommand{\fp}[2]{{}_{\hspace{3pt}#1\hspace{-2pt}}\times_{#2}}

\DeclareMathOperator{\id}{id}
\DeclareMathOperator{\expl}{Expl}
\DeclareMathOperator{\coker}{coker}

\providecommand{\et}[2]{\ensuremath{\bold T^{#1}_{#2}}}
\providecommand{\lrb}[1]{\ensuremath{\left(#1\right)}}
\providecommand{\abs}[1]{\left\lvert #1\right\rvert}

 \thanks{This paper was written during the authorÕs stay at the Mathematical Science Research Institute in Berkeley.} 
\title{De Rham theory of exploded manifolds}  

\begin{document}

\begin{abstract}

This paper extends de Rham theory of smooth manifolds to exploded manifolds. Included are versions of Stokes' theorem, de Rham cohomology, Poincar\'e duality, and integration along the fiber. The resulting de Rham cohomology theory of exploded manifolds is used in a separate paper \cite{egw} to define Gromov--Witten invariants of exploded manifolds.

\end{abstract}
\maketitle

\section{Introduction}

The goal of this paper is to describe a version of de Rham cohomology for exploded manifolds which extends de Rham cohomology for smooth manifolds. At first glance the most natural extension of de Rham cohomology would be to take the complex of smooth or $\Cs\infty1$ differential forms on an exploded manifold with the usual differential $d$. Unfortunately, this naive extension does not have good properties - for example, in a smooth connected family of exploded manifolds, the cohomology defined this way might change. Moreover, the tools of integration and Poincar\'e duality are not available for this naive extension. 

Instead, we shall use a sub-complex $\Omega^{*}(\ex B)$ of $\Cs\infty1$ differential forms on $\ex B$, defined below in definition \ref{omega def}. In the case that $\ex B$ is a smooth manifold, $\Omega^{*}(\ex B)$ is the usual complex of smooth differential forms. We shall show in Section \ref{invariant section} that the cohomology $H^{*}(\ex B)$ does not change in connected families of exploded manifolds. (This fact is nontrivial to prove because families of exploded manifolds are not always locally trivial.) As suggested by the names of the sections of this paper, many of the standard tools of de Rham cohomology still apply for $\Omega^{*}(\ex B)$.

\

From now on, some knowledge of the definitions and notation from \cite{iec} shall be necessary to understand this paper. 
Recall that coordinates on $\mathbb R^{n}\times \et mP$ are given by \[x_{j}\co \mathbb R^{n}\times \et mP\longrightarrow \mathbb R\text{ for }1\leq j\leq n\] and \[\tilde z_{i}\co \et mP\longrightarrow \mathbb C^{*}\e{\mathbb R}\text{ for }1\leq i\leq m\] Smooth or $\Cs\infty1$ differential  one-forms on $\mathbb R^{n}\times\et mP$ are given by smooth or $\Cs\infty1$ functions times $dx_{j}$ and the real and imaginary parts of $\tilde z_{i}^{-1}d\tilde z_{i}$. These differential forms are not ideal  for de Rham cohomology, as even compactly supported forms may not have finite integral.

\begin{example}[A compactly supported form with infinite integral]\end{example}

 $\et 11 :=\et 1{[0,\infty)}$ has a single coordinate \[\tilde z\co \et 11\longrightarrow \mathbb C^{*}\e{[0,\infty)}\] 

 Consider the two-form $\alpha$ given by the wedge product of the real and imaginary parts of $\tilde z^{-1} d\tilde z$. Over any tropical point $\e a\in \totb{\et 11}$ in the tropical part of $\et 11$, there is a $\mathbb C^{*}$ worth of points corresponding to a choice of coefficient $c$ of $\tilde z=c\e a$. On the $\mathbb C^{*}$ worth of points over each tropical point of $\totb{\et 11}$, $\alpha$ is a nonzero $\mathbb C^{*}$ invariant volume-form, so by any straightforward definition of integration, $\alpha$ should have infinite integral. Similarly, if $\alpha$ is multiplied by any continuous function $f\co \et 11\longrightarrow \mathbb R$ which is nonzero when $\totl{\tilde z}=0$, the integral of $f\alpha$ is again infinite. This is because $f\alpha$ restricted to the $\mathbb C^{*}$ worth of points over any point in $(0,\infty)\subset\totb{\et 11}$ is a nonzero $\mathbb C^{*}$--invariant volume-form, and hence has infinite integral. 
 
 Recall that
 \[\totl{c\e a}=\begin{cases}0\text{ if }a>0
 \\ c\text{ if }a=0
 \end{cases}\]
 and that the topology on $\et 11$ is a non-Hausdorff topology in which every open subset is the pullback of some open subset in $\mathbb C$ under the map $\totl{\tilde z}\co \et 11\longrightarrow \mathbb C$. It follows that $f\alpha$ may be compactly supported and still have infinite integral.

\

 There are several possible fixes to this problem - we shall consider forms which do not contain the real part of $\tilde z_{i}^{-1}d\tilde z_{i}$ where it is an obstacle to integration. In particular, we shall require that our differential forms vanish on integral vectors, which are the vectors $v$ so that $vf$ is an integer times $f$ for all exploded functions $f$. (For example, the integral vectors on $\et 11$ are integers times the real part of $2\tilde z\frac \partial{\partial\tilde z}$ wherever $\totl{\tilde z}=0$.)

\

For Stokes' theorem to work out correctly, we shall also require the following condition: Given any map  $f\co \et 1{(0,\infty)}\longrightarrow \ex B$, we shall assume that our differential forms vanish on all of  the vectors in the image of $Tf$.

As an example to see that some restriction is necessary for Stokes' theorem to hold, consider a compactly supported form $
\theta$ on $\et 11$ given by the imaginary part of $\tilde z^{-1}d\tilde z$ multiplied by a smooth, compactly supported  function $f$ which is $1$ when $\totl{\tilde z}=0$. Then the integral of $d\theta$ over $\et 11$ is $2\pi$ rather than $0$.

\begin{defn}\label{omega def}Let  $\Omega^{k}(\ex B)$ be the vector-space of $\Cs\infty1$ differential $k$--forms $\theta$ on $\ex B$ so that
\begin{enumerate}
\item\label{od1} for all integral vectors $v$, the differential form $\theta$ vanishes on $v$,
\item\label{od2} and for all maps $f\co \et 1{(0,\infty)}\longrightarrow \ex B$, the differential form $\theta$ vanishes on all vectors in the image of $Tf\co T\et 1{(0,\infty)}\longrightarrow T\ex B$. 
\end{enumerate}

Denote by $\Omega^{k}_{c}(\ex B)\subset \Omega^{k}(\ex B)$ the subspace of forms with complete support. (Say that a form has complete support if  the set where it is non zero is contained inside a complete subset of $\ex B$ --- in other words, a compact subset with tropical part consisting only of complete polytopes.)

\end{defn}

Clearly, the usual wedge product, exterior differential, and interior product with a $\Cs\infty 1$ vectorfield are all defined and obey the usual properties on $\Omega^{*}(\ex B)$. Moreover, given any $\Cs\infty1$ map $g\co \ex B\longrightarrow \ex C$, the pullback $g^{*}$ of differential forms sends forms in  $\Omega^{*}(\ex C)$ to forms in $\Omega^{*}(\ex B)$. This is because  $Tg$ always sends integral vectors to integral vectors, and sends any vector in the image of $Tf\co T\et 1{(0,\infty)}\longrightarrow T\ex B$ to a vector in the image of $T(g\circ f)$.

\begin{defn} Denote the homology of $(\Omega^{*}(\ex B),d)$ by $H^{*}(\ex B)$, and the homology of $(\Omega^{*}_{c}(\ex B),d)$ by $H^{*}_{c}(\ex B)$.
\end{defn}

We shall show in section \ref{poincare section} that given an assumption about the topology of $\ex B$ akin to the existence of a finite good cover,  $H^{*}_{c}(\ex B)$ is dual to $H^{*}(\ex B)$.

\section{Mayer Vietoris sequence}

 Below we shall prove that the usual Mayer Vietoris sequence holds. This requires partitions of unity, which are constructed in section \ref{partition}. 

\begin{lemma}[Mayer Vietoris sequence]\label{Mayer Vietoris} Given open subsets $U$ and $V$ of an exploded manifold $\ex B$, the Mayer Vietoris sequences
\[0\longrightarrow \Omega^{*}(U\cup V)\xrightarrow{\theta\mapsto \theta\oplus \theta} \Omega^{*}(U)\oplus \Omega^{*}(V)\xrightarrow{\theta_{1}\oplus\theta_{2}\mapsto\theta_{1}-\theta_{2}} \Omega^{*}(U\cap V)\longrightarrow 0 \]
\[0\longrightarrow \Omega^{*}_{c}(U\cap V)\xrightarrow{\theta\mapsto\theta\oplus-\theta} \Omega^{*}_{c}(U)\oplus\Omega^{*}_{c}(V)\xrightarrow{\theta_{1}\oplus\theta_{2}\mapsto\theta_{1}+\theta_{2}} \Omega^{*}_{c}(U\cup V)\longrightarrow 0\]
are exact sequence of chain complexes. 
\end{lemma}

\pf 

The proof is identical to the proof in the case of smooth manifolds given in \cite{BT}. We shall discuss the first exact sequence first.

As usual in the Mayer Vietoris sequence, the first map is the direct sum of the restriction of forms from $U\cup V$ to $U$ and $V$, which is an injective chain map. Then the second map is the restriction of forms on $U$ to $U\cap V$ minus the restriction of forms from $V$ to $U\cap V$. This is a chain map, and its kernel is the forms which agree on $U\cap V$, which obviously agrees with the image of the first map. It remains to verify that this second map is surjective. Choose a partition of unity for $U\cup V$ subordinate to   $U$ and $V$, so we have smooth functions $\rho_{U}$ and $\rho_{V}$ on $U\cup V$ which sum to $1$ and which are supported inside $U$ and $V$ respectively. Then any  form $\theta\in \Omega^{*}(U\cap V)$ is in the image of $\rho_{V}\theta\oplus(-\rho_{U}\theta)\in \Omega^{*}(U)\oplus\Omega^{*}(V)$.

Now for the second exact sequence. The first map is given by inclusion of completely  supported forms in $U\cap V$ to completely supported forms in  $U$ and $V$. This is clearly an injective chain map. The second map is given by inclusion of completely supported forms in $U$ to $U\cup V$ plus the inclusion of completely supported forms in $V$ to $U\cup V$. Again, it is clear that this is a chain map. The kernel consists of forms which cancel each other on $U\cap V$, and which are also supported in $U\cap V$. This agrees with the image of the first map. To see that the second map is surjective, suppose that $\theta\in \Omega^{*}(U\cup V)$. Then $\theta$ is the image of $\rho_{U}\theta\oplus\rho_{V}
\theta\in\Omega^{*}(U)\oplus\Omega^{*}(V)$. 

\stop

\section{Integration and Stokes' theorem}

We shall show below that if $\ex B$ is oriented and $n$--dimensional, then the integral of compactly supported forms in $\Omega^{*}(\ex B)$ is well defined. 

\

\begin{example}[$\Omega^{2}(\et 11)$] \end{example}
On $\et 11$, the integral vectors are integer multiples of twice the real part of $\tilde z\frac\partial{\partial\tilde z}$ at points where $\totl{\tilde z}=0$. Any two-form in $\Omega^{2}(\et 11)$ must therefore vanish wherever $\totl{\tilde z}=0$.  What remains is the subset of $\et 11$ where $\tilde z\in\mathbb C^{*}$. Let $\totl{\tilde z}=e^{r+i\theta}$ and  denote the imaginary part of $\tilde z^{-1}d\tilde z$  by $d\theta$ and denote the real part by $dr$. If $\alpha\in\Omega^{2}(\et 11)$, then
\[\alpha=f(r,\theta)dr\wedge d\theta\]
where $f$ is a smooth function of $r$ and $\theta$, and for any $\delta<1$, the size of $f$ or any of  its derivatives is bounded by $e^{\delta r}$ as $r\to -\infty$. The form $\alpha$ is compactly supported on $\et 11$ if and only if $f$ vanishes when $r$ is sufficiently large. (Of course, $\alpha$ need not be compactly supported in $\mathbb C^{*}\subset \et 11$ to  be compactly supported in $\et 11$.) The integral of $\alpha$ is finite if $\alpha$ is compactly supported in $\et 11$ and given by 
\[\int_{\et 11}\alpha=\int_{-\infty}^{\infty}\int_{0}^{2\pi}f(r,\theta)drd\theta\ .\]

\

Any top-dimensional form in $\Omega^{*}(\ex B)$ will vanish on all strata apart from those strata of $\ex B$ which are smooth manifolds (and therefore have no nonzero integral vectors). We can therefore define the integral of a top-dimensional form $\theta$ on an oriented exploded manifold $\ex B$ to be the sum of the integrals of $\theta$ over these smooth strata. This integral is well defined if the integral over each smooth stratum is well defined and the resulting sum of integrals is well defined.

\begin{defn}
If $\alpha$ is a top-dimensional form on an oriented exploded manifold $\ex B$, define the integral of $\alpha$ to be the sum of the integral of $\alpha$ over all strata of $\ex B$ which are smooth manifolds. 
\[\int_{\ex B}\alpha=\sum_{\totl{\ex B_{i}} =\text{ point}}\int_{\ex B_{i}}\alpha\]
 
\end{defn}

\begin{lemma} \label{finite integral} If  a top-dimensional form $\alpha\in\Omega^{*}(\ex B)$ is compactly supported, then the integral of $\alpha$ is finite. 
\end{lemma}

\pf

By using a partition of unity, we may assume that $\alpha$ is compactly supported within a single coordinate chart $\mathbb R^{n}\times \et mP$. The smooth manifold strata of this coordinate chart  are the strata over the (zero-dimensional) corners of the polytope $P$. As $P$ has only finitely many such corners, we need only verify the finiteness of our integral over one stratum of our coordinate chart. 

  We must deal with the possiblity that our corner of $P$ may not be standard. Pulling back $\alpha$ to a refinement of $\mathbb R^{n}\times \et mP$ will not change the integral. We may subdivide our corner of $P$ so that the corresponding corner of each new cell has exactly $m$ edges. Therefore, by taking a refinement and again using a partition of unity, we may assume that a neighborhood of the corner of $P$ at our stratum is  isomorphic to a neighborhood of $0$ in the image of some integral-affine map applied to  the standard quadrant $[0,\infty)^{m}$. It follows that our stratum is contained in the image of a proper map from $\mathbb R^{n}\times \et m{[0,\epsilon)^{m}}$ to our coordinate chart, and that this map restricted to our stratum is a covering map of some positive degree. 

It therefore suffices to prove our lemma for a compactly supported form  $\alpha\in\Omega^{n+2m}(\mathbb R^{n}\times \et m{[0,\epsilon)^{m}})$.  Use coordinates $\totl{\tilde z_{k}}=e^{r_{k}+i\theta_{k}}$ on our stratum. Then
\[\alpha=f(x,r,\theta)\prod dx_{j}\prod dr_{k}\wedge d\theta_{k}\]
where $f$ is smooth and bounded by some constant times $e^{\frac 12\sum_{k}r_{k}}$. Furthermore,
 on the support of $f$, $\abs x$ and $r$ are bounded above.    The integral of $\alpha$ on our stratum is therefore finite and well defined.

\stop

\

Define an exploded manifold with boundary to be an abstract exploded space $\ex M$ locally isomorphic to $(-\infty,0]\times\mathbb R^{n}\times \et mP$. As usual, if $\ex M$ is oriented, the boundary $\partial \ex M$ is oriented in a way consistent with giving the boundary of  $(-\infty,0]\times\mathbb R^{n}\times \et mP$ the usual orientation on $\mathbb R^{n}\times \et mP$, so a positively oriented top-dimensional form on $\partial \ex M$ is obtained by inserting an outward pointing normal vector into a top-dimensional positively oriented form on $\ex M$. We can now state Stokes' theorem for exploded manifolds:

\begin{thm}[Stokes' theorem] \label{stokes} If $\ex M$ is an oriented exploded manifold with boundary and $\theta\in\Omega^{*}_{c}(\ex M)$, then 
\[\int_{\ex M}d\theta=\int_{\partial \ex M}\theta\]
\end{thm}

\pf

We shall use Stokes' theorem for smooth manifolds. Because of the linearity of the equation we must prove, we may use a partition of unity to reduce to the case when $\ex M$ is covered by a single coordinate chart. Consider the integral of $d\theta$ over a single stratum $M'$ of our coordinate chart. We must deal with the following problem: even though $d\theta$ is compactly supported on $\ex M$, it may not have compact support on $ M'$.

 The tropical part of $ M'$ is a $0$--dimensional corner of the polytope $\totb{\ex M}$. Identify $\totb{\ex M}$ with a polytope  $P\subset [0,\infty)^{m}$ so that $\totb{ M'}$ is $0$. Using the corresponding map $\ex M\longrightarrow \et m{[0,\infty)^{m}}$ we may consider the coordinates $\tilde z_{i}$ from $\et m{[0,\infty)^{m}}$ as coordinates on $\ex M$. Consider the hypersurface 
$N_{\epsilon}\subset M'$ where $\abs{\tilde z_{1}\dotsc\tilde z_{m}}=\epsilon$, oriented as the boundary of the region $M_{\epsilon}'$ where   $\abs{\tilde z_{1}\dotsc\tilde z_{m}}\geq\epsilon$. Our form $\theta$ is compactly supported when restricted to  $M'_{\epsilon}$, so we can use Stokes' theorem for manifolds.  
\[\int_{M'}d\theta=\lim_{\epsilon\to 0}\int_{M'_{\epsilon}}d\theta=\int_{\partial M'}\theta+\lim_{\epsilon\rightarrow 0}\int_{N_{\epsilon}}\theta  \]

We must consider  the integral $\int_{N_{\epsilon}}\theta$ as $\epsilon\to 0$. More radically,   consider setting $\epsilon =1\e x$ for $x>0$ small enough that the hypersurface $ x_{1}+\dotsb+x_{m}=x$ intersects each stratum of $P$ attached to $0$. In this case $N_{1\e x}:=\{\abs{\tilde z_{1}\dotsc\tilde z_{m}}=1\e x\}\subset \ex M$ may be regarded as an infinite union of hypersurfaces, one over each point in $P$ where $\sum_{i}x_{i}=x$. We shall argue below that, of these hypersurfaces,  only the hypersurfaces $N_{e}$ over the edges $e$ of $P$ shall contribute to $\int_{N_{\epsilon}}\theta$ as $\epsilon\to 0$. (With the tool of integration over the fiber, (Theorem \ref{fiber integration}), it becomes clear that $\int_{N_{\epsilon}}\theta$ depends continuously on $\epsilon$ where $\epsilon\in \mathbb R^{*}\e{\mathbb R}$ is given the topology induced from some refinement of $\et 1{[0,\infty)}$. This is because we can first integrate $\theta$ along the fiber of the map $\tilde z_{1}\dotsb\tilde z_{n}\co \ex M\longrightarrow \et 1{[0,\infty)}$, then observe that the integral around a circle $\abs{\tilde z}=\epsilon$ of a form on a refinement of $\et 1{[0,\infty)}$  depends continuously on $\epsilon$ where $\epsilon\in \mathbb R^{*}\e{\mathbb R}$, and $\mathbb R^{*}\e{\mathbb R}$ is given the topology from our refinement of $\et 1{[0,\infty)}$. In what follows, we give a more basic argument of this fact.)

Because $\theta $ is compactly supported, it vanishes when any smooth monomial is large enough.
Because $\theta$ vanishes on integral vectors and is one less than top-dimensional, $\theta$ vanishes on all strata of $\totb{\ex M}$ with tropical dimension at least $2$. Therefore, $\theta$ is bounded by a constant times $w_{\mathcal S}^{\delta}$ for any $\delta<1$ where $\mathcal S$ consists of all strata of $\totb{\ex M}$ with dimension at least $2$. (Recall from \cite{iec} that $w_{\mathcal S}$ is a finite sum of absolute values of smooth monomials which vanish on all the strata in $\mathcal S$.) Note that $w_{\mathcal S}$ is a finite sum of exponentials, (and is not constant along any straight line within $N_{\epsilon}$) so the integral of $\theta$ over the regions where $w_{\mathcal S}$ is small will also be small. 
It follows that in the limit $\epsilon\to 0$ the integral $\int_{N_{\epsilon}}\theta$ is concentrated in the directions corresponding to the edges of $\totb{\ex M}$ attached to our stratum, and that 

\[\lim_{\epsilon\rightarrow 0}\int_{N_{\epsilon}}\theta=\sum_{\text{edges } e}\int_{N_{e}}\theta \]
where $N_{e}$ is   the hypersurface where $\abs{\tilde z_{1}\dotsc\tilde z_{m}}=1\e x$  in the stratum corresponding to an edge $e$ attached to our corner. (The sum is over all these edges.) As with $N_{\epsilon}$,  $N_{e}$ is oriented as the boundary of the region where $\abs{\tilde z_{1}\dotsc\tilde z_{m}}\geq 1\e x$. The integral of $d\theta$ over our coordinate chart will therefore have a contribution for each end of each edge of $\totb {\ex M}$. Consider the two contributions corresponding to an edge $e$ of $\totb {\ex M}$ which has two ends. These contributions will be the integral of $\theta$ over two different hypersurfaces in the stratum corresponding to $e$. Each of these hypersurfaces is in the form of $\{\abs {\tilde z^{\alpha} }=c\}$, and is transverse to the integral vectors in this stratum. Because $\theta$ vanishes on integral vectors, and is constant in the direction of integral vectors, the map given by projection along the direction of the integral vectors between these two  hypersurfaces  preserves $\theta$, but reverses the orientation of these two hypersurfaces. Therefore, the contributions from the two different ends of $e$ cancel each other out. The facts that the support of $\theta$ is complete and that $\theta$ vanishes on the image of any $\et 1{(0,\infty)}$, imply that if $\theta$ does not vanish on $N_{e}$, there must be a corner of $\totb{\ex M}$ at both ends of the edge $e$. Therefore all the contributions from $\int_{N_{e}}\theta$ cancel, and we obtain that
\[\int_{\ex M}d\theta=\int_{\partial M}\theta\]
as required.

\stop

\section{Cohomology of a coordinate chart}

In this section, we calculate $H^{*}$ for all standard coordinate charts and $H^{*}_{c}$ for coordinate charts $\mathbb R^{n}\times \et mP$ for which $P$ is complete.

\begin{lemma}\label{poincare} Let $P\subset \mathbb R^{m}$ be an integral-affine polytope. Suppose that the directions of the infinite rays in $P$ span the last $k$ coordinate directions in $\mathbb R^{m}$. Then  $H^{*}(\mathbb R^{n}\times \et mP)$ is equal to the free exterior algebra generated by the imaginary parts of $\tilde z_{i}^{-1}d\tilde z_{i}$ for $i=1,\dotsc, m-k$.
\end{lemma}

\pf

Note that if $M$ is a smooth manifold, $H^{*}(M)$ is the usual De Rham cohomology of $M$.
There is an obvious map of a torus $\mathbb T^{m}\longrightarrow \mathbb R^{n}\times \et mP$  which pulls back the above forms to non trivial homology classes in $H^{*}(\mathbb T^{m})$, so $H^{*}(\mathbb R^{n}\times \et mP)$ must contain a copy of the free exterior algebra generated by the above differential forms.

We shall prove that each class in $H^*(\mathbb R^{n}\times \et mP)$ may be represented by a differential form which is constant in standard coordinates and contains no $dx_{i}$ factors. Our proof shall follow the proof of the Poincar\'e lemma in \cite{BT}.

Consider the map $K\co \Omega^{*}\longrightarrow \Omega^{*-1}$ given by 
\[K(\theta)(x_{1},x_{2},\dotsc )=\int_{0}^{x_{1}}i_{\frac \partial{\partial x_{1}}}\theta (s,x_{2},\dotsc)ds\] 
then 
\[dK(\theta)=(1-dx_{1}\wedge i_{\frac \partial{\partial x_{1}}})\lrb{\int_{0}^{x_{1}} di_{\frac\partial{\partial x_{1}}}\theta(s,x_{2}\dotsc) ds}+ dx_{1}\wedge i_{\frac\partial{\partial x_{1}}}\theta\]
so
\[(Kd+dK)(\theta)=\int_{0}^{x_{1}}(i_{\frac\partial{\partial x_{1}}}d+(1-dx_{1}\wedge i_{\frac \partial{\partial x_{1}}})di_{\frac\partial{\partial x_{1}}})\theta(s,x_{2}\dotsc) ds+ dx_{1}\wedge i_{\frac\partial{\partial x_{1}}}\theta\]
Suppose that $i_{\frac \partial{\partial x_{1}}}\theta=0$, then 
\[(Kd+dK)(\theta)=\int_{0}^{x_{1}}L_{\frac\partial{\partial x_{1}}}\theta(s,x_{2},\dotsc)ds=\theta-\theta(0,x_{2},\dotsc)\]
Suppose that $\theta=dx_{1}\wedge\alpha$ where $i_{\frac\partial{\partial x_{1}}}\alpha=0$. Then 
\[(Kd+dK)dx_{1}\wedge\alpha=\int_{0}^{x_{1}}0ds+dx_{1}\wedge\alpha\]
Therefore, in general

\[(Kd+dK)\theta=\theta-(1-dx_{1}\wedge i_{\frac\partial {\partial x_{1}}})\theta(0,x_{2},\dotsc,)\]

It follows that we can represent the cohomology class of any closed form $\theta$ with the closed form $(1-dx_{1}\wedge i_{\frac\partial {\partial x_{1}}})\theta(0,x_{2},\dotsc,)$ which is independent of $x_{1}$ and $dx_{1}$. Similarly, we may represent any class in $H^{*}(\mathbb R^{n}\times\et mP)$ by a form pulled back from $\et mP$ under the obvious projection map.

Now we have reduced to the case of differential forms on $\et mP$. The standard basis for  differential forms on $\et mP$ is given by the exterior algebra generated by the real and imaginary parts of $\tilde z_{i}^{-1}d\tilde z_{i}$, so we can consider forms in $\Omega^{*}(\et mP)$ as maps from $\et mP$ to $\mathbb R^{2^{2m}}$. (Of course, not all $\Cs\infty1$ maps to $\mathbb R^{2^{2m}}$ will correspond to forms in $\Omega^{*}$ because of the condition that forms in $\Omega^{*}$ vanish on integral vectors.) We wish to show that any cohomology class can be represented by a form which is constant in this basis. In particular, if $P^{\circ}$ is the interior stratum of $P$, then we shall show that a closed form $\theta$ represents the same cohomology class as $e_{P^{\circ}}\theta$. (Here we use notation from the definition of $\Cs\infty1$ from section $7$ of \cite{iec}.)

For each stratum of $P$, choose an integral vector $\alpha$ pointing towards the interior of $P$. Choose these vectors consistently so that the vectors for adjacent strata differ only by a vector contained within one of the strata. One way to do this is to rescale and deform $P$ until all its vertices are on integer points within $\mathbb R^{n}$ and $0$ is in its interior. Then for each stratum $S$  of $P$, we may choose $\alpha$ to be the negative of some integral point which is contained in the closure of $S$.

\

\begin{figure}
\includegraphics{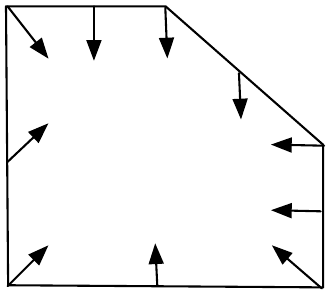}
\end{figure}

\

On each stratum, consider the vector field 
\begin{equation}\label{v def}v:=\sum_{i=1}^{m} \alpha_{i}\frac \partial {\partial r_{i}}\end{equation}
where $\frac\partial {\partial r_{i}}$ is the real part of $2\tilde z_{i}\frac \partial{\partial z_{i}}$. This vector field $v$ is not a globally defined vectorfield on $\et mP$ because it may change by an integral vector from one stratum of $\et mP$ to the next. We can think of $v$ as a `vectorfield defined up to integral vectors'. As differential forms $\theta\in \Omega^{*}(\et mP)$ always vanish on integral vectors, $i_{v}\theta$ is still a well defined form in $\Omega^{*}(\et mP)$ even though $v$ might jump by an integral vectorfield when changing from stratum to stratum. Let $\Phi_{tv}$ be the flow of the vectorfield $v$  for time $t$ on each stratum. We will not be overly worried by the fact that $\Phi_{tv}$ does not give a globally defined map from $\et mP$ to itself. Note that the flow of any integral vectorfield preserves any $\Cs\infty1$ differential form, therefore the ambiguity in the definition of $v$ does not affect how forms are changed by the flow of $v$. Therefore, if $\theta\in \Omega^{*}(\et mP)$, then  $\Phi_{tv}^{*}\theta\in\Omega^{*}(\et mP)$ for all $t$. 

Given $\theta\in\Omega^{*}(\et mP)$ define
\begin{equation}\label{defK}K\theta:=\int_{-\infty}^{0}\Phi_{tv}^{*}i_{v}\theta dt\end{equation} 

 Now  check that $K\theta\in\Omega^{*}$. Note that $i_{v}\theta$ vanishes on $\et m{P^{\circ}}\subset \et mP$ because on this stratum $v$ is an integral vector field. Also note that $\Phi_{tv}$ travels towards this central stratum as $t\rightarrow -\infty$. In particular, given any smooth monomial  $\zeta$ on $\et mP$, $L_{v}\abs{\zeta}$ is some positive number times $\abs \zeta$ . It follows that $\Phi_{tv}^{*}\abs \zeta$ restricted to any compact subset of a stratum is bounded by some constant times $e^{t}$ as $t\rightarrow -\infty$.   The fact that $i_{v}\theta$ vanishes on the central stratum and is $\Cs\infty1$ implies that  restricted to any compact subset, there exists a constant $c$ so that 
\begin{equation}\label{Kestimate}\abs{\Phi_{tv}^{*}i_{v}\theta} <ce^{\frac 12 t} \end{equation}
and similar estimates hold for any derivative of $\theta$. Therefore, on each stratum $K\theta$ is well defined and smooth. It is clear that $K\theta$ vanishes on all vectors that differential forms in $\Omega^{*}$ should vanish on, so it remains to check that $K\theta$ is $\Cs\infty1$. At this stage, the reader must be familiar with section 7 of  \cite{iec}.

Note that for any stratum $S$, 
$e_{S}\Phi_{tv}^{*}i_{v}\theta=\Phi_{tv}^{*}i_{v}e_{S}\theta$. Therefore for any collection of strata $\mathcal S$, $\Delta_{S}K\theta=K\Delta_{S}\theta$. We must prove that $w_{\mathcal S}^{-\delta}\Delta_{\mathcal S}K\theta$ is bounded on any compact subset for any $\delta<1$ (and we must prove a similar estimate for any derivative of $\theta$).  For any $0<\epsilon<1-\delta$, $\Delta_{\mathcal S}\theta$ is bounded by a constant times $w_{S}^{\delta+\epsilon}$ on any compact subset. The weight $w_{\mathcal S}$ is a sum of absolute values of smooth monomials which vanish on $P^{\circ}$, so $\Phi_{tv}^{*}w_{\mathcal S}$ is bounded on compact subsets by a constant times $e^{t}$ for $t<0$. It follows that $\Phi^{*}_{tv}\Delta_{S}i_{v}\theta$ is bounded by a constant times $w^{\delta}_{\mathcal S}e^{\epsilon t}$ on compact subsets for $t<0$. Integrating this gives that  $w_{\mathcal S}^{-\delta}\Delta_{\mathcal S}K\theta$ is bounded on any compact subset. For any constant vectorfield $X$, note that $L_{X}K\theta=KL_{X}\theta$, so the bounds for the derivatives of $\theta$ follow from the same argument, and $K\theta\in\Omega^{*}(\et mP)$. 

Now consider
\begin{equation}\label{Kprop}(Kd+dK)\theta=\int_{-\infty}^{0}\Phi^{*}_{tv}(i_{v}d+di_{v})\theta dt=\theta-\lim_{t\rightarrow -\infty}\Phi^{*}_{tv}\theta=\theta-e_{P^{\circ}}\theta\end{equation}
 It follows that if $\theta\in\Omega^{*}(\et mP)$ is any closed differential form, we may represent the same cohomology class by the constant differential form $e_{P^{\circ}}\theta$. The lemma follows as the only constant differential forms in $\Omega^{*}(\et mP)$ are the forms mentioned in the statement of the lemma.
 
 \stop 

\

Note that $(\et 11)^{m}$ has the same cohomology as $\mathbb C^{m}$. Using a good cover, (constructed in Lemma \ref{good cover}), we may use this to prove that the cohomology of the explosion of any compact complex manifold relative to a normal crossing divisor is equal to the cohomology of the original manifold.

\begin{cor}
If $M$ is a compact complex manifold with a normal crossing divisor, then the smooth part map
\[\expl M\xrightarrow {\totl\cdot} M\]
induces an isomorphism on cohomology.
\[\totl\cdot^{*}\co H^{*}(M)\xrightarrow{\cong} H^{*}(\expl M)\]
More generally, if $\ex B$ has a finite good cover in the sense of Lemma \ref{good cover}, and all polytopes in the tropical part of $P$ are quadrants $[0,\infty)^{m}$ then $\totl\cdot^{*}$ is an isomorphism on cohomology.
\[\totl\cdot^{*}\co H^{*}(\totl{\ex B})\xrightarrow{\cong} H^{*}(\ex B)\]

\end{cor}

\pf 

Note that the map 
\[\totl{\cdot}\co \expl M\longrightarrow M\]
may be considered as a smooth map of exploded manifolds, where $M$ is given the structure of a smooth manifold, regarded as an exploded manifold. In this setting $\Omega^{*}M$ just corresponds to the usual smooth differential forms on the smooth manifold $M$, and we may pullback differential forms in $\Omega^{*}$ and cohomology classes as usual for exploded manifolds.

Choose a finite good cover $\{U_{i}\}$ of $\ex B$ in the sense of Lemma \ref{good cover}, so the intersection of any number of these $U_{i}$ is either empty or isomorphic to $\mathbb R^{n}\times \et mP$. By assumption the  only possible polytopes $P$ are quadrants $[0,\infty)^{m}$.  Lemma \ref{poincare} tells us that $H^{*}(\mathbb R^{n}\times (\et 11)^{m})$ is generated by the constant functions. The smooth part of $\mathbb R^{n}\times (\et 11)^{m}$ is $\mathbb R^{n}\times \mathbb C^{m}$, so the smooth part map $\totl\cdot$ induces an isomorphism on cohomology:

\[\totl\cdot^{*}\co H^{*}(\mathbb R^{n}\times \mathbb C^{m})\xrightarrow{\cong} H^{*}(\mathbb R^{n}\times (\et 11)^{m})\]

Therefore, if $\ex B$ has a good cover by a single open set, our lemma holds. We may now proceed by induction over the cardinality of a good cover using the Mayer Vietoris sequence from Lemma \ref{Mayer Vietoris}. 

Suppose that our lemma holds for all exploded manifolds satisfying our tropical part assumption with a good cover containing at most $k$ sets. Then suppose that $\ex B$ has a good cover $\{U,V_{1},\dotsc,V_{k}\}$. Let $V=\bigcup_{i=1}^{k} V_{k}$. Then our lemma holds for $U$, $V$ and $U\cap V$.  Then the smooth part map gives the following commutative diagram involving Mayer Vietoris sequences

{\setlength\mathsurround{0pt}
\[\begin{tikzcd}\Omega^{*}\lrb{U\cup V}\rar&\Omega^{*}(U)\oplus \Omega^{*}(V)\rar& \Omega^{*}(U\cap V)
\\ \uar{\totl\cdot^{*}}\Omega^{*}\lrb{\totl{U}\cup \totl V}\rar &\uar{\totl\cdot^{*}\oplus \totl\cdot^{*}}\Omega^{*}(\totl U)\oplus \Omega^{*}(\totl V)\rar & \Omega^{*}(\totl U\cap \totl V)\uar{\totl{\cdot}^{*}}
\end{tikzcd}\]
}

Considering the induced maps on the homology long exact sequence and using the five lemma then implies that 
\[\totl\cdot^{*}\co H^{*}(\totl U\cup \totl V)\longrightarrow H^{*}(U\cup V)\] 
is an isomorphism. By induction, our lemma must hold for  $\ex B$ so long as $\ex B$ has a finite good cover and the tropical part of $\ex B$ contains only quadrants. The tropical part of the $\expl M$ contains only quadrants, and Lemma \ref{good cover} implies that if $M$ is compact, $\expl M$ has a finite good cover, so our lemma also holds for $\expl M$.

\stop 

\

It is not true in general that $\totl{\ex B}$ has the same cohomology as $\ex B$. For example, $\et 1{[0,1]}$ has the same cohomology as $\mathbb C^{*}$, but the smooth part of $\et 1{[0,1]}$ is two copies of $\mathbb C$ glued at $0$.

\

Stokes' theorem implies that   $(\alpha,\theta)\mapsto\int_{\ex B}\alpha\wedge \theta$ gives a bilinear pairing $H^{k}_{c}(\ex B)\times H^{n-k}(\ex B)\longrightarrow \mathbb R$ where the dimension of $\ex B$ is $n$. We shall prove that in many situations, this pairing is non-degenerate, so Poincar\'e duality holds. For this, we will need an assumption on the polytopes $P$ contained in the tropical part of $\ex B$ --- we will need that these polytopes are complete, contain no entire lines, and are simplicial at infinity; see Definition \ref{simplicial}. Rougly speaking, this assumption is equivalent to $P$ being a deformation of a cone spanned by linearly independent vectors; it is also satisfied by any compact polytope.

\

Lemma \ref{local poincare} below computes $H^{*}_{c}(\mathbb R^{n}\times \et mP)$ for complete polytopes $P$ which are simplicial at infinity and contain no entire lines. This computation also implies the computation of $H^{*}_{c}$ in the case where $P$ is a complete polytope which is simplicial at infinity, because there exists an obvious projection map $\pi\co \mathbb R^{n}\times \et mP\longrightarrow \mathbb R^{n}\times \et {m'}{P'}$ so that $P'$ is complete, simplicial at infinity,  and contains no lines, and so that that $\pi^{*}$ is a bijection on both $\Omega^{*}$ and $\Omega^{*}_{c}$. 

\begin{lemma}\label{local poincare}Suppose that $P$ is a complete polytope which is simplicial at infinity and contains no entire lines. Then the integration pairing
\[(\alpha,\theta)\mapsto\int_{\mathbb R^{n}\times \et mP}\alpha\wedge\theta\] identifies  $H^{*}_{c}(\mathbb R^{n}\times\et mP)$ with the dual of $H^{*}(\mathbb R^{n}\times \et mP)$.

\end{lemma}

\pf

 By averaging, we may represent any class in $H^{*}_{c}$ by a differential form $\theta$  preserved by the flow of the imaginary part of $\tilde z_{i}\frac \partial {\partial \tilde z_{i}}$. As $P$ is simplicial at infinity, there exist $k$ linearly-independent integral vectorfields $b'_1,\dotsc, b'_k$ on $P$ whose nonnegative span contains all integral vectorfields whose flow exists for all positive time. Each such vectorfield $b_i'=(b'_{i,1},\dotsc,b'_{i,m})$ on $P$ corresponds to a vectorfield $b_i$ on $\mathbb R^{n}\times\et mP$ that is the imaginary part of $\sum_j b_{i,j}\tilde z_{j}\frac \partial {\partial \tilde z_{j}}$. These vectorfields $b_i$ have the following significance: a form $\theta$ satisfying Definition \ref{omega def} part \ref{od1} satisfies part \ref{od2} if and only if $i_{b_i}\theta=0$ above the strata of $P$ that $b'$ is tangent to. As these $b_i$ are linearly independent, we can choose a basis $\beta_1,\dots,\beta_m$ for the span of the imaginary parts of $\tilde z_{i}^{-1}d\tilde z_{i}$ such that $i_{b_i}\beta_j$ is $1$ if $i=j$ and is otherwise zero. (These $\beta_i$ correspond to constant one-forms $\beta_i'$ on $P$.) With this choice of basis, each torus-invariant closed differential form in $\Omega^*(\mathbb R^{n}\times \et mP)$ breaks up into a sum of closed forms $\alpha\wedge \beta\in \Omega^*(\mathbb R^{n}\times \et mP)$ where $\alpha$ vanishes on the imaginary part of $\tilde z_i\frac\partial{\partial \tilde z_i}$, and $\beta$ is a product of these $\beta_j$. Below, we will define a map $f\co \mathbb R^{n}\times \et mP\longrightarrow \mathbb R^{n+m}$ and replace $\alpha$ with a standard form pulled back from $\mathbb R^{n+m}$.

Choose a basis $\{\zeta_{v}\co =\totl{\e{a}\tilde z^{v}}\}$ for the smooth monomials on $\et mP$.
 Recall that $\zeta_{v}$ is a smooth $\mathbb C$--valued function so $\abs{\zeta_{v}}$ is a smooth positive function. 
Consider the differential form
\[d\left(\sum_{i=1}^{n}\frac 12\abs{ x_{i}}^{2}+\sum_{v}\abs{\zeta_{v}}^{2}\right)\] 
as giving a smooth map 
\[f\co \mathbb R^{n}\times \et mP\longrightarrow \mathbb R^{n+m}\]
with first $n$ components the insertion of $\frac\partial {\partial x_{i}}$, and the last $m$ components the insertion of the real part of $\tilde z_{i}\frac\partial{\partial \tilde z_{i}}$. We shall check below that that this map $f$ is proper,  and the image of a stratum $S$ with $k$--dimensional tropical part  is a cone $C_{S}$ in $\mathbb R^{n+m}$ with codimension $k$. (Throughout this proof, we shall use $S$ to refer to  a stratum of $\mathbb R^{n}\times \et mP$, and $\totb{S}$ to refer to  its tropical part which is a stratum of the polytope $P$.)
\begin{figure}
\includegraphics{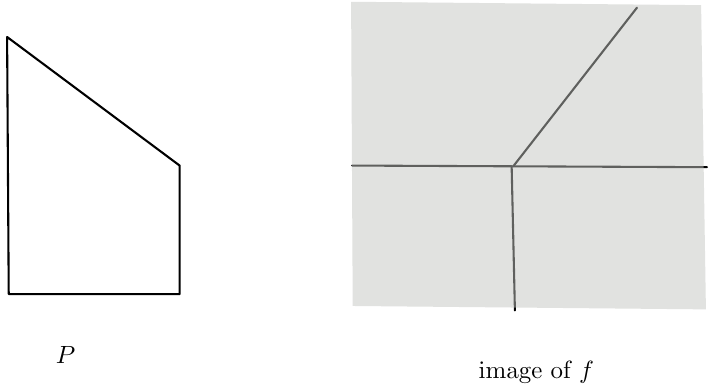}
\end{figure}

Let us examine the map $f$. 
\[f=\id_{\mathbb R^{n}}+\sum_{v}\abs{\zeta_{v}}^{2}v\]
In the above, identify $v\in\mathbb Z^{m}$ with the corresponding vector in $0\times \mathbb R^{m}\subset\mathbb R^{n}\times \mathbb R^{m}$. This formula implies that $f$ is proper and that the image of $f$ restricted to a stratum $S$ is contained inside the cone, $C_{S}$, defined as  $\mathbb R^{n}$ times the positive span of all $v$ so that $\zeta_{v}$ is nonzero on our stratum.
\[C_{S}\co =\mathbb R^{n}\times \text{Positive span}\{v \text{ so that $\zeta_{v}\neq 0$ on $S$}\}\subset \mathbb R^{n+m}\]
 Taking the derivative of $f$ gives
\[Df=\id_{\mathbb R^{n}}+\sum_{v}2\abs{\zeta_{v}}^{2}\abs{v}^{2}\pi_{v}\] 
where \[\abs v^{2}\pi_{v}=(\sum_{i=1}^{m} v_{i}dr_{i})v\]
and $dr_{i}$ indicates the real part of $\tilde z_{i}^{-1}d\tilde z_{i}$. If we regard the vector-space with basis $\{\frac \partial {\partial x_{i}},\frac \partial {\partial r_{j}}\}$ as $\mathbb R^{n+m}$, then $\pi_{v}$ indicates orthogonal projection onto the subspace spanned by $v$. Regarded this way, $Df$ is a symmetric, positive semidefinite matrix. It follows from this formula that restricted to any stratum $S$, $Df$ is surjective onto the tangent space to the cone $C_{S}$. Combined with the properness of $f$, this implies that the image of $f$ restricted to $S$ is the interior of this cone $C_{S}$. 

So, the image of $f$ is $\mathbb R^n$ times the dual fan of $P\subset\mathbb R^m$, with the image of $S$ being  $\mathbb R^n$ times the cone dual to $\totb S\subset P$ --- in other words, the set of linear functions $P\longrightarrow \mathbb R$ which are constant on $\totb S$, and are strictly greater than this everywhere off the closure of $\totb S$. Each $\beta_i$ corresponds to a vector $\beta_i'$ in the image of this dual fan, and the image of $f$ is $\mathbb R^n$ times the non-negative span of $\beta'_1,\dotsc,\beta'_k$ times the span of $\beta'_{k+1},\dotsc, \beta'_m$. So, $b_1',\dotsc,b'_k$ define linearly independent functions on $\mathbb R^{n+m}$ such that the image of $f$ is the set where $b'_i\geq 0$. Morever, for $i\leq k$, and a smooth form $\theta$ on $\mathbb R^{n+m}$,  the form $f^*\theta\wedge\beta_i $ is in $\Omega^*(\mathbb R^{n}\times \et mP)$ if and only if  $\theta$ vanishes where $b'_i=0$.

Restricted to a stratum $S$ with $k$--dimensional tropical part $\totb S$, the formula for $Df$ implies that the fiber of $f\co S\longrightarrow \mathbb R^{n+m}$ over any point in $C_{S}\subset \mathbb R^{n+m}$ is equal to $\et k{\totb S}$ times a $(m-k)$--dimensional torus and has tangent space spanned by integral vectors and the imaginary parts of the $\tilde z_{i}\frac\partial{\partial \tilde z_{i}}$.   Our differential form $\alpha$ vanishes restricted the tangent space to these fibers, and is closed, therefore $\alpha$ restricted to any stratum $S$ must be equal to the pullback under $f$ of some differential form, $\alpha'$, that is smooth on the interior of $C_{S}$. In general, it will not be true that $\alpha'$ comes from a smooth differential form on $\mathbb R^{n+m}$; for that, we must modify $\alpha$.

Consider the operator $K$ defined in equation (\ref{defK}) on page \pageref{defK}. As noted in equation (\ref{Kprop}), this oporator has the property that $d(K\alpha)=\alpha-e_{P^{\circ}}\alpha$. Choose some compactly supported smooth function $\rho$ on $\et mP$ which is $1$ in a neighborhood of the interior stratum $P^{\circ}$ of $\et mP$, and modify $\alpha$ to the form $\alpha-d(\rho K\alpha)$. This modified form (which we shall again call $\alpha$) is still compactly supported, but has the property that in a neighborhood of this interior stratum, $\alpha=e_{P^{\circ}}\alpha$ is the pullback of some smooth form via the composition of the map $f$ with the orthogonal projection to $C_{P^{\circ}}$. 

 Suppose that for all strata $\totb S$ of $P$ with dimension greater than $k$, there exists a neighborhood of $S$ on which $\alpha$ is the pullback of a smooth form under the composition of $f$ with the orthogonal projection to $C_{S}$. We shall now modify $\alpha$ so that the same holds for the strata $\totb S$ with dimension $k$. Let $F$ indicate the smallest face of $P$ which contains $\totb S$ (in other words $F$ is the closure of $\totb S\subset P$). Using the implicit function theorem for exploded manifolds proved in \cite{iec}, we may identify a neighborhood of our  stratum $S$  with $\mathbb R^{n}\times (\mathbb C^{*})^{m-k}\times \et {k}F$ so that $f$ composed with orthogonal projection to $C_{S}$ is equal to $e_{S}f$ and so that these neighborhoods for different strata $\totb S$ of dimension $k$ do not intersect. Let $K$ be the operator defined in equation (\ref{defK}) for these new coordinates. To define $K$ we use the flow of a vector field $v$ from equation (\ref{v def}); the flow of the vector field $v$ commutes with $e_{S}$, so   our inductive hypothesis  ensures   that the form $K\alpha$ vanishes on the intersection of the open set where it is defined with a neighborhood of the strata with higher dimensional tropical part, because on this region,
  $i_{v}\alpha=0$. Let $\rho$ be a compactly supported function on $\et kF$ which is $1$ in a neighborhood of $\et k{\totb S}\subset \et kF$. Then $\rho K\alpha$ is a compactly supported form, and we may modify $\alpha$ 
to $\alpha-d(\rho K\alpha)$ without changing its cohomology class in $H^{*}_{c}$. Doing the same for all strata $\totb S$ of dimension $k$, we get a modified $\alpha$ which satisfies the  required condition because on a neighborhood of $S$,  $e_{S}\alpha=\alpha$, and because $\alpha$ has not been modified  on a neighborhood of all the strata of higher tropical dimension.

We may therefore modify $\alpha$ so that each stratum $S$ has a neighborhood so that $\alpha$ is the pullback of a smooth form under $f$ composed with orthogonal projection to $C_{S}$. It follows that this modified form $\alpha$ is $f^{*}\alpha'$ for some smooth closed form $\alpha'$ on $\mathbb R^{n+m}$. 

 We now have that $H^{*}_{c}$ is generated by differential forms $f^{*}\alpha'\wedge \beta$ where $\beta$ is some product of $\beta_{i}$. Choose some standard form $\alpha_{0}$ with integral $1$ which is compactly supported in the interior of the image of $f$. We shall now show that we may exchange $\alpha'$ for $(\int{\alpha'})\alpha_{0}$. Recall that the $\beta_i$ for $i\leq k$ are our basis forms not vanishing on the image of some $\et 1{(0,\infty)}$.  Because $\alpha'\wedge \beta$ must vanish on all vectors in the image of any $\et 1{(0,\infty)}$, if $\beta$ contains $\beta_1\dotsb \beta_k$, then $\alpha'$ must vanish on the boundary of $\bigcup_{S}C_{S}$. (Because $P$ is complete, every point in the boundary of $\bigcup_{S}C_{S}$ is the image of $f$ composed with some nontrivial map from $\et 1{(0,\infty)}$, and $\beta_1\dotsb\beta_k$ never vanishes on all vectors in the image of such a map.)  As $\alpha'$ in a neighborhood of such a boundary is the pullback of some form via projection to the boundary, it follows that  $\alpha'$ is compactly supported inside the interior of $\bigcup_{S}C_{S}$. As  $\bigcup _{S}C_{S}$ is the dual fan to $P$ which is a closed polytope in $[0,\infty)^{n+m}$, the interior of $\bigcup_{S}C_{S}$  diffeomorphic to $\mathbb R^{n+m}$. Therefore \[\alpha'-d\gamma=\alpha_{0}\int\alpha'\] where $\gamma$ is compactly supported inside the interior of $\bigcup_{S}C_{S}$. As $f^{*}\gamma\wedge \beta\in \Omega^{*}_{c}$, the modified form $(\int\alpha')f^{*}\alpha_{0}\wedge\beta$ represents the same class in $H^{*}_{c}$ as $\alpha\wedge\beta$. 

Now consider the case of $f^*\alpha'\wedge \beta$ where $\beta$ does not contain $\beta_1\dotsb\beta_k$. Recall that the image of $f$ is $\mathbb R^n$ times the positive span of $\beta'_1,\dotsc,\beta'_k$ times the span of $\beta'_{k+1},\dotsb \beta'_m$, or, the subset where the linear functions $b'_i\geq 0$ for $i=1,\dotsc,k$. If $\beta$ contains $\beta_i$ for $i\in I\subsetneq \{1,\dotsc,k\}$, but does not contain $\beta_j$ for $j\notin I$,  then $\alpha'$ is reqired to vanish where $\prod_{i\in I}b'_i=0$, but not required to vanish on the rest of the boundary of the image of $f$. There is a homotopy of $\mathbb R^{n+m}$ in the direction of $-\beta'_j$ so that the support of $\alpha'$ is moved outside the image of $f$. This homotopy can be used to construct  a  compactly supported form $\gamma$ on $\mathbb R^{n+m}$ which vanishes on $\prod_{i\in I}b'_i$  and for which $d\gamma =\alpha'$ on the image of $f$. As $f^{*}\gamma\wedge \beta\in\Omega^{*}_{c}$, we have that  $\alpha\wedge\beta$ represents the zero cohomology class in $H^{*}_{c}$.

In conclusion, we have shown that $H^{*}_{c}$ is generated by forms 
\[f^{*}\alpha_{0}\wedge \beta_1\dotsb \beta_{k}\wedge\beta\]
 where $\beta$ is some product of $\beta_{i}$ for $k<i\leq m$. Lemma \ref{poincare} showed that $H^{*}$ is generated as an exterior algebra by $\beta_{i}$ for $k<i\leq m$. The integration pairing on our space of differential forms times $H^{*}$  is therefore nondegenerate, therefore all the above forms represent independent classes in $H^{*}_{c}$, and the integration pairing identifies $H^{*}_{c}$ with the dual of $H^{*}$ as required.

\stop

\section{Poincar\'e duality}\label{poincare section}

 \begin{thm}[Poincar\'e duality]\label{Poincare duality}
 If $\ex B$ is a complete oriented exploded manifold that is simplicial at infinity\footnote{Definition \ref{simplicial}.},  so that each map $\ex T\longrightarrow \ex B$ is constant, then the integration pairing  gives an isomorphism  between $H^{*}(\ex B)$ and its dual. 
 
 \
 
 More generally, if $\ex B$  has a finite good cover in the sense of Lemma \ref{good cover} and each polytope in the tropical part of $\ex B$ is complete, simplicial at infinity,  and contains no entire lines, then the integration pairing identifies  $H^{*}(\ex B)$ with the dual of $H^*_{c}(\ex B)$.
 
 \end{thm}
 
 \pf
 
 We shall use Lemma \ref{good cover}, which states that any complete exploded manifold $\ex B$ must have a finite good cover by open sets so that any intersection is isomorphic to a standard coordinate chart $\mathbb R^{n}\times \et mP$. The condition that each map $\ex T\longrightarrow \ex B$ is constant implies that the polytope $P$ contains no entire lines, so we may apply  Lemma \ref{local poincare} to know that the integration pairing identifies the dual of $H^{*}(\mathbb R^{n}\times \et mP)$ with $H^{*}_{c}(\mathbb R^{n}\times \et mP)$. The proof may now proceed as in the case of smooth manifolds by induction over the size of our open cover using the Mayer Vietoris sequences from Lemma \ref{Mayer Vietoris}. 
 
 In particular, suppose that the dimension of $\ex B$ is $n$. Define the differential $d'=(-1)^{n+1-k}d$ on $\Omega_{c}^{k}$. The Mayer Vietoris sequence obviously is still exact for this new differential, and the homology of $(\Omega_{c}^{*},d')$ is obviously the same as the homology of $(\Omega_{c}^{*},d)$. This sign modification allows the following formula:
 \[\int_{\ex B}(d\alpha)\wedge \beta=\int_{\ex B}\alpha\wedge d'\beta\]
Let $C_{*}$ denote the dual chain complex to $(\Omega_{c}^{*},d')$. The above formula implies that the integration pairing on any oriented $n$--dimensional manifold gives a chain map $\Omega^{*}\longrightarrow C_{*}$. We shall now verify that the corresponding map between Mayer Vietoris sequences is commutative
\[\begin{array}{ccccccccc}0&\longrightarrow&\Omega^{*}(U\cup V)&\longrightarrow &\Omega^{*}(U)\oplus \Omega^{*}(V)&\longrightarrow &\Omega^{*}(U\cap V)&\longrightarrow&  0
\\ &&\downarrow &&\downarrow &&\downarrow
\\0&\longrightarrow&C_{*}(U\cup V)&\longrightarrow &C_{*}(U)\oplus C_{*}(V)&\longrightarrow &C_{*}(U\cap V)&\longrightarrow&  0 \end{array}\]
 
 Following $\alpha\in \Omega^{*}(U\cup V)$ across and down, then evaluating on $\beta_{1}\oplus \beta_{2}\in \Omega^{*}_{c}(U)\oplus\Omega^{*}_{c}(V)$ gives 
 \[\int_{U}\alpha\wedge \beta_{1}+\int_{V}\alpha\wedge \beta_{2}\]
 Following $\alpha$ down and across, then evaluating on $\beta_{1}\oplus\beta_{2}$ gives
 \[\int_{U\cup V}\alpha\wedge(\beta_{1}+\beta_{2})\]
 Therefore the first square commutes. The commutativity of the second square  amounts to the equation
 \[\int_{U\cap V}(\alpha_{1}-\alpha_{2})\wedge\beta=\int_{U}\alpha_{1}\wedge \beta+\int_{V}\alpha_{2}\wedge(-\beta)\]
  where $\beta$ is compactly supported in $U\cap V$.
  
  Therefore, taking homology gives a commutative diagram 
  \[\begin{array}{ccccccccc}H^{k-1}(U\coprod V)&\hspace{-6mm}\rightarrow\hspace{-6mm} &H^{k-1}(U\cap V)&\hspace{-6mm}\rightarrow\hspace{-6mm}&H^{k}(U\cup V)&\hspace{-6mm}\rightarrow\hspace{-6mm} &H^{k}(U\coprod V)&\hspace{-6mm}\rightarrow\hspace{-6mm} &H^{k}(U\cap V)
 \\ \downarrow & & \downarrow & & \downarrow & & \downarrow &&\downarrow
 \\ \hspace{-3mm}H^{n+1-k}_{c}(U\coprod V)^{*}\hspace{-3mm}&\hspace{-6mm}\rightarrow\hspace{-6mm} &\hspace{-1mm}H^{n+1-k}_{c}(U\cap V)^{*}\hspace{-3mm}&\hspace{-6mm}\rightarrow\hspace{-6mm}&\hspace{-1mm} H^{n-k}_{c}(U\cup V)^{*}\hspace{-3mm}&\hspace{-6mm}\rightarrow\hspace{-6mm} &\hspace{-1mm}H^{n-k}_{c}(U\coprod V)^{*}\hspace{-3mm}&\hspace{-6mm}\rightarrow\hspace{-6mm} &\hspace{-1mm}H_{c}^{n-k}(U\cap V)^{*}\end{array}\]
  The downward arrows above are given by the integration pairing. Say that Poincar\'e duality holds if this integration pairing map is an isomorphism. The Five Lemma implies that if Poincar\'e duality holds on $U$, $V$ and $U\cap V$, then Poincar\'e duality holds on $U\cup V$. 
 
 Suppose that Poincar\'e duality holds on all oriented exploded manifolds satisfying our assumptions on their tropical part and having a good cover containing at most $k$ members. Suppose that $\ex B$ satisfies the tropical part assumptions and has a good cover $\{U_{1},\dotsc,U_{k+1}\}$. Then Poincar\'e duality must hold for $U_{k+1}$, $\bigcup_{i=1}^{k}U_{i}$ and $U_{k+1}\cap\bigcup_{i=1}^{k}U_{i}$. The above argument then gives that Poincar\'e duality must hold for $\ex B$. By induction starting with Lemma \ref{local poincare}, Poincar\'e duality must hold for all oriented exploded manifolds that admit  
a finite good cover and which have a tropical part containing only complete polytopes that contain no lines. Lemma \ref{good cover} states that complete exploded manifolds have a finite good cover, and our theorem follows.
  
 \stop 
 
 \

Poincar\'e duality as stated in Theorem
\ref{Poincare duality} does not imply the usual relationship between intersections of submanifolds and wedge products of Poincar\'e duals. 
We shall explore this relationship further when we return to Poincar\'e duality in section \ref{pdfp}.

 \section{Integration along the fiber}
 
 In this section, we define integration along the fibers, $f_{!}$, for suitable maps $f$. Given a $\Cs\infty1$ map, $f\co \ex A\longrightarrow \ex B$, and a compactly supported differential form $\theta$ on $\ex A$, we may regard $\theta$ as a current (something dual to the space of differential forms), then push forward this current to obtain a current on $\ex B$. If this current on $\ex B$ is also represented by a differential form, we call this differential form $f_{!}(\theta)$. In particular, when it exists, $f_{!}\theta$ has the  property that for all differential forms $\alpha$ on $\ex B$, 
 
 \[\int_{\ex B}\alpha\wedge f_{!}\theta=\int_{\ex A}(f^{*}\alpha)\wedge \theta\]
  
  In the case of oriented smooth manifolds, $f_{!}$ exists if $f$ is a submersion. In our case, we must be careful that $Tf$ restricted to the subspace spanned by integral vectors is also surjective.
  
  \begin{thm}\label{fiber integration} Let $\ex A$ and $\ex B$ be oriented exploded manifolds, and   suppose that  $f\co \ex A\longrightarrow \ex B$ is a $\Cs\infty1$ map satisfying  the following conditions:
  \begin{enumerate}
  \item f is a submersion in the sense that
  \[Tf\co T_{x}\ex A\longrightarrow T_{f(x)}\ex B\] is surjective.
\item \[Tf\co {}^{\mathbb Z}T_{x}\ex A\longrightarrow {}^{\mathbb Z}T_{(f(x))}\ex B\] is a surjective map on integral vectors.

\end{enumerate}
Then, if the fiber of $f$ is $n$--dimensional, there exists a linear chain map $f_{!}\co  \Omega^{*}_{c}(\ex A)\longrightarrow \Omega^{*-n}_{c}(\ex B)$, with the property that 
\[\int_{\ex B} \alpha\wedge f_{!}\theta=\int_{\ex A}f^{*}\alpha\wedge \theta\]
for all $\alpha\in\Omega^{*}(\ex B)$. 

If all polytopes $P$ in the tropical part of $\ex B$ are complete and contain no entire lines, then $f_{!}(\theta)$ is uniquely determined by the above property.

  \end{thm}
 
  \pf
  
   The discussion on fiber-products in section 9 of \cite{iec} implies that each fiber of $f$ is a $\Cs\infty1$ exploded manifold and that the top wedge of the cotangent space of the fibers is a $\Cs\infty1$ vectorbundle, $\bigwedge^{\text{top}}T^{*}_{\text{vert}}\ex A$, over $\ex A$. The pullback of $\bigwedge T^{*}\ex B$ is also a vectorbundle, $f^{*}\bigwedge T^{*}\ex B$, over $\ex A$, and may be regarded as a sub bundle of $\bigwedge T^{*}\ex A$ because $f$ is a submersion.
   The tensor of these two bundles is a $\Cs\infty1$ vectorbundle \[E :=\bigwedge^{\text{top}}T^{*}_{\text{vert}}\ex A\otimes f^{*}\bigwedge T^{*}\ex B \] over $\ex A$. From $\theta\in\Omega^{*}_{c}(\ex A)$, we can associate a $\Cs\infty1$ section $\theta'$ of $E$ as follows: inserting a top-dimensional polyvector $v$ tangent to the fiber of $f$ into the righthand places of $\theta$  gives a form $\theta\wedge v$ which vanishes on the kernel of $Tf$. Therefore $\theta\wedge v$ must be a section $f^{*}\bigwedge T^{*}\ex B$ considered as a sub bundle of $\bigwedge T^*\ex A$. Define $\theta'$ as the unique section of $E$ so that $\theta'(v)=\theta\wedge v$. It is obvious that this definition of $\theta'$ does not actually depend on the choice of top-dimensional polyvector $v$. As $\theta$ has complete support and vanishes on all the vectors which forms in $\Omega^{*}$ must vanish on, $\theta'$ restricted to the fiber $f^{-1}(p)$ is in $\Omega^{n}_{c}(f^{-1}(p))\otimes \bigwedge T^{*}_{p}\ex B$. Orient the fibers of $f$ so that if $\alpha$ is a volume-form on $\ex B$ and $\beta$ is a volume-form on the fibers of $f$, then $f^{*}\alpha\wedge \beta $ is a volume-form on $\ex A$.
  
  We may therefore integrate $\theta'$ along the fiber $f^{-1}(p)$ to obtain a form $f_{!}\theta(p)\in \bigwedge T^{*}_{p}\ex B$. We must now verify that $f_{!}\theta$ defined this way is in $\Omega^{*}_{c}(\ex B)$, and verify that it satisfies our defining property.
  
  As $f$ is a submersion, any $\Cs\infty1$ vector field $v$ on $\ex B$ lifts to a $\Cs\infty1$ vectorfield $\tilde v$ on $\ex A$ so that $Tf(\tilde v)=v$. Let $\Phi_{t\tilde v}$ indicate the flow of the vectorfield $\tilde v$ on $\ex A$ and $\Phi_{tv}$ indicate the flow of the vectorfield $v$ on $\ex B$. We have that $f\circ\Phi_{t\tilde v}=\Phi_{tv}\circ f$ and $f_{!}\circ \Phi_{t\tilde v}^{*}=\Phi_{tv}^{*}\circ f_{!}$. As the map $f_{!}$ is linear, and $\Phi_{t\tilde v}^{*}(\theta)$ is differentiable in $t$, $f_{!}\Phi_{t\tilde v}^{*}\theta$ is also differentiable in $t$ and   
   \[f_{!}L_{\tilde v}(\theta)=L_{v}f_{!}\theta\ .\]

Given a vector field $v$ on $\ex B$, note that the section $(i_{\tilde v}\theta)'$ of our bundle $E$ does not depend on the choice of lift $\tilde v$. If at $p$, $v$ is an integral vectorfield, then our second assumption on $f$ implies that given any point $q\in f^{-1}(p)$,  we may choose our lift $\tilde v$ so that $\tilde v$ is an integral vectorfield at $q$, therefore $(i_{\tilde v}\theta)'$ vanishes around $q$. Therefore $(i_{\tilde v}\theta)'$ vanishes on $f^{-1}(p)$, so $f_{!}\theta$ vanishes on integral vectors.  Similarly, given any map of $\et 1{(0,\infty)}$ passing through $p=f(q)$, the fact that $\theta$ has complete support and $Tf$ is surjective on integral vectors  implies that either $\theta(q)=0$ or this map may be covered by a map $g$ of $\et 1{(0,\infty)}$ to $\ex A$ composed with $f$ so that the image of $g$ contains $q$. It follows that if $v$ is in the image of such a map, $(i_{\tilde v}\theta)'$ must vanish. Therefore $f_{!}\theta$ must vanish on the tangent space of the image of any map from $\et 1{(0,\infty)}$. Therefore, $f_{!}\theta$ vanishes on all vectors it should vanish on.

As the image of any complete set is complete, $f_{!}\theta$ has complete support. 
Therefore, to check that $f_{!}\theta\in\Omega^{*}_{c}{\ex B}$, it remains to check that $f_{!}\theta$ is $\Cs\infty1$.  To do this, we work locally in a single coordinate chart $U'$ on $\ex A$ and $U$ on $\ex B$. Our assumptions on $f$ imply that the image of every stratum of $\totb{U'}$ under $\totb f$ is a stratum of $\totb{U}$. 

Recall, from section 7 of \cite{iec}, that if $\zeta$ is a monomial function on a coordinate chart and $S$ is a stratum, then $e_{S}\zeta$ is $0$ if $\zeta$ vanishes on the stratum $S$, and otherwise $e_{S}\zeta=\zeta$. A function $g$ on our coordinate chart may be regarded as a function of smooth monomials $\zeta_{i}$, and $e_{S}g(\zeta_{1},\dotsc,\zeta_{m})=g(e_{S}\zeta_{1},\dotsc,e_{S}\zeta_{m})$. The reader should also keep in mind that $e_{S}g$ at a point $p$ in our coordinate chart should be regarded as sampling $g$ at a point $q$ shifted from $p$ to be in the stratum $S$. (To make sense of $e_{S}\theta$ where $\theta$ is a differential form, the standard basis for differential forms in a coordinate chart is used in order to identify differential forms at different points.)

By modifying our chart $U'$ using the implicit function theorem if necessary, we may assume that the pullback of monomial functions from $U$ are monomial functions on $U'$. It follows that given any stratum $S$ of $U'$ and $\Cs\infty1$ function $g$ on $\ex B$, $e_{S}(f\circ g)=f\circ e_{\totb f(S)}g$.   Note that $f_{!}\theta$ depends only on position in $\totl{U}$; if $p$ and $p'$ have the same image in $\totl{U}$, then the integrals used to compute $f_{!}{\theta}$ are the same on the fiber over $p$ and $p'$. It follows that $e_{S}f_{!}\theta$ makes sense.

Given any set $\mathcal S$ of strata of $\totb{U}$, let $\mathcal S':=\totb f^{-1}(\mathcal S)$ be the set of strata of $\totb{U'}$ sent to $\mathcal S$ by $\totb f$. Consider the case of a single stratum, $S$, of $U$. Recall that $e_{S}f_{!}\theta$ at a point $p$ is equal to $e_{S}f_{!}\theta$ at a point $q$ obtained from $p$ by shifting $p$ into the stratum $S$.  The fiber of $f$ over $q$ intersects all the strata $T$ in $ S':=\totb f^{-1}(S)$. For any stratum $T$ whose intersection with $f^{-1}(q)$ is a manifold,  $e_{T}\theta'$ restricted to $f^{-1}(p)$ may be regarded as the pullback of $\theta'$ under a diffeomorphism of $f^{-1}(p)$ with $T\cap f^{-1}(q)$.  Accordingly, 
\[e_{S}f_{!}\theta=\sum_{T\in S'}f_{!}(e_{T}\theta)\ .\]

  If $T_{1}$ and $T_{2}$ are distinct strata in $S'$, then $e_{T_{1}}e_{T_{2}}\theta'=0$ because it samples $\theta'$ on a stratum on which there are integral vectors in the vertical tangent bundle. Therefore 
  \[\Delta_{S}f_{!}\theta:=(1-e_{S})f_{!}\theta=f_{!}\lrb{\prod_{T\in S'}(1-e_{T})\theta}:= f_{!}(\Delta_{S'}\theta)\ .\]
It follows that for any set $\mathcal S$ of strata of $\totb U$, 
\[\Delta_{\mathcal S}f_{!}\theta= f_{!}(\Delta_{\mathcal S'}\theta)\ .\]  
Recall from \cite{iec} that the definition of a $\Cs\infty1$ function $g$ involves controlling $w_{\mathcal S}^{-\delta}\Delta_{S}D^{n}g$ for all $0<\delta<1$, sets of strata $\mathcal S$ and number of derivatives $n\geq 0$. We need to be able to compare the weighing function $w_{\mathcal S}$ on $U$ with the corresponding weighting function $w_{\mathcal S'}$ on $U'$. Recall that $w_{\mathcal S}$ is defined as a sum of absolute values of monomial functions $\zeta$ that vanish on all strata in $\mathcal S$. Each monomial function $\zeta$ is the smooth part of some exploded  function in the form $\e x\tilde z^{\alpha}$; the monomial function $\zeta$ vanishing on $S$ is equivalent to the integral-affine function $x+\alpha$ on $\totb U$ being nonnegative on $\totb U$, and strictly positive on $S$.
  
We may restrict to the case that the closure of the support of $\theta$ intersects the interior stratum of $U'$. As $\theta$ has complete support, this implies that every nonnegative integral-affine function, $\alpha$ on $\totb {U'}$, achieves its minimum on every fiber of $\totb{U'}\longrightarrow \totb U$. Choose a point $y\in \totb{U'}$ so that $\alpha$ achieves its fiberwise minimum at $y$, and $y$ projects to the interior of $\totb U$.   Let $F$ be the smallest face of $\totb{U'}$ containing $y$. Then  $\alpha$  is fiberwise constant on $F$ and achieves its  fiberwise minimum on $F$. The  assumption that $Tf$ is surjective on integer vectors implies that $\alpha$ restricted to $F$ is equal to the pullback of some nonnegative integral-affine function from $\totb {U}$.  It follows that any nonnegative integral-affine function on $\totb{U'}$ which is positive on all strata in $\mathcal S'$ is equal to the pullback of some integral-affine function which is positive on $\mathcal S$, plus some nonnegative integral-affine function. It follows that every smooth monomial which vanishes on all strata in $\mathcal S'$ is divisible by the pullback of a smooth monomial vanishing on all strata in $\mathcal S$. Therefore, we may choose $w_{\mathcal S'}=f^{*}w_{\mathcal S}$. It follows that 
  \[w_{\mathcal S}^{-\delta}\Delta_{S}f_{!}\theta=f_{!}(w_{\mathcal S'}^{-\delta}\Delta_{\mathcal S'}\theta)\]
  As we already have that $L_{v}f_{!}\theta=f_{!}L_{\tilde v}\theta$, it follows that $f_{!}\theta$ is $\Cs\infty1$ if $\theta\in \Omega^{*}_{c}$. So $f_{!}\theta\in \Omega^{*}_{c}(\ex B)$ if $\theta\in \Omega^{*}_{c}(\ex A)$. 
  
We have defined our map $f_{!}$ in the same way as the integration over fibers map for smooth manifolds with the sign convention $\int \alpha\wedge f_{!}\theta=\int f^{*}\alpha\wedge \theta$. (See for example \cite{BT}.) As our integrals are just defined as a finite sum of integrals over smooth manifolds and this formula holds for smooth manifolds, it also holds for us:

\[\int_{\ex A}f^{*}\alpha\wedge \theta=\int_{\ex B}\alpha\wedge f_{!}\theta\] 

The above formula uniquely characterizes $f_{!}$ in the case that $\ex A$ and $\ex B$ are smooth manifolds. A quick calculation using Stokes theorem gives that 
\[\int_{\ex A}f^{*}\alpha\wedge d\theta=\int_{\ex B}\alpha \wedge d f_{!}\theta\] 
therefore, in the case of smooth manifolds $f_{!}$ is a chain map. As our $f_{!}$ is simply obtained by a sum of $f_{!}$ for smooth manifold components, our $f_{!}$ is also a chain map.
  
  \stop
  
  \begin{remark}The map $f_{!}\co \Omega_{c}^{*}\ex A\longrightarrow \Omega_{c}^{*}\ex B$ only depends on the relative orientation of $\ex A$ and $\ex B$; if the opposite orientations of $\ex A$ and $\ex B$ are used, then $f_{!}$ remains unchanged. We can define $f_{!}$ in the case that we have a choice of relative orientation for $f$  so that $f_{!}$ locally coincides with the map from Theorem \ref{fiber integration} when we locally choose an orientation of $\ex B$. See also Remark \ref{relative orientation}
  \end{remark}
  
  \section{fiber-products and integration along the fiber}
  
  In Lemma \ref{fp integration} below, we shall show that $f_{!}$ transforms well under fiber-products. 
  In order to do this, we need to specify the orientation convention we shall use for fiber-products. Fiber-products of exploded manifolds are defined in \cite{iec}. 
 It is also shown in \cite{iec} that if $f$ and $g$ are transverse, then  the derivatives of the maps in the following commutative diagram

{\setlength\mathsurround{0pt}
\[\begin{tikzcd}\ex A\fp fg\ex B\rar{f'}\dar{g'} &\ex B\dar{g}
    \\ \ex A\rar{f} &\ex C \end{tikzcd}\]}
give a short exact sequence
\[0\longrightarrow T_{(p_{1},p_{2})}\lrb{\ex A\fp fg\ex B}\xrightarrow{(Tg',Tf')} T_{p_{1}}\ex A\times T_{p_{2}}\ex B\xrightarrow{Tf-Tg} T_{f(p_{1})}\ex C\longrightarrow 0\]
In other words, we have the same relationship between tangent spaces as in the case of manifolds, so we may orient fiber-products of exploded manifolds as we orient fiber-products of manifolds.  In particular,  the above exact sequence and commutative diagram imply that 
\begin{itemize}
\item $Tf'$ gives an isomorphism between $\ker Tg'$ and $\ker Tg$
 \item $Tf$ gives  an isomorphism between $\coker Tg'$ and $\coker Tg$
 \item $Tg'$ gives an isomorphism between $\ker Tf'$ and $\ker Tf$
 \item $Tg$ gives an isomorphism between $\coker Tf'$ and $\coker Tf$
 \item $\ker T(g\circ f')=\ker Tf'\oplus\ker Tg' \cong \ker Tf\oplus\ker Tg
$ 
\item $\coker T(g\circ f')\cong \coker Tf\oplus \coker Tg$. 
\end{itemize}

Use the convention that if $V$ and $W$ are oriented vector-spaces so that $\sigma_{V}$ and $\sigma_{W}$ are positive forms, their direct sum $V\oplus W$ is oriented so that $\sigma_{V}\wedge \sigma_{W}$ is positive. With this convention, a choice of orientation of any two of $V$, $W$ and $V\oplus W$ implies an orientation on the third.

An orientation of $V$ relative to $W$ is an choice of orientation of $V\oplus W$; this should be regarded as giving an orientation of $V$ for any choice of orientation of  $W$. For example, a choice of isomorphism $\phi\co V\longrightarrow W$ gives a natural orientation of $V$ relative to $W$ so that $\phi$ is oriented.

  Given a map $A\co X\longrightarrow Y$, between oriented vector-spaces, there are several possible conventions for orienting $\ker A$ relative to $\coker A$. The following definition gives notation for describing some of these conventions. 

\begin{defn}
Given a map of oriented vector-spaces $A\co X\longrightarrow Y$ we shall use the following shorthand for an orientation convention for $\ker A$ relative $\coker A$. By saying the identification
\[\coker A\oplus X=\ker A\oplus Y\]
 is an oriented isomorphism, we mean that given any metric on $X$ and $Y$, the natural map
 \[A'\co \coker A\oplus X\longrightarrow \ker A\oplus Y\]
is an oriented isomorphism. This map $A'$  is defined to restrict to $\coker A$ to be the identification of $\coker A$ with the orthogonal complement of $A(X)\subset Y$, and  restrict to $X$ to be the  orthogonal projection onto $\ker A$ and the map $A$.\end{defn}

Of course, the relative orientation of $\ker A$ and $\coker A$ given by  the isomorphism $A'$ does not depend on the choice of metrics on $X$ and $Y$.

Taking the above direct sums in different orders gives different orientation conventions.

\

We shall find the following way of arranging kernels and cokernels convenient.

\[\coker Tf \oplus T(\ex A\fp fg \ex B) \oplus \coker Tg= \ker Tf\oplus T\ex C\oplus \ker Tg \] 
 
  \begin{defn}[Orientation convention for fiber-products] 
  
  Let $\ex A$, $\ex B$ and $\ex C$ be oriented exploded manifolds, and let $f\co \ex A\longrightarrow \ex C$ and $g\co \ex B\longrightarrow \ex C$ be transverse maps.
  Orient  $\ker Tf$ relative to $\coker Tf$ so that the identification
  \[\coker Tf \oplus T\ex A=\ker Tf\oplus T\ex C\]
  is an oriented isomorphism. On the other hand, orient $\ker Tg$ relative to $\coker Tg$ so that the following identification gives an oriented isomorphism:
  \[T\ex B\oplus \coker Tg=T\ex C\oplus \ker Tg\]
  
  Then orient $T\ex A\fp fg \ex B$ so that the following identification is an oriented isomorphism:
  \[\coker Tf \oplus T(\ex A\fp fg \ex B) \oplus \coker Tg= \ker Tf\oplus T\ex C\oplus \ker Tg \] 

  \end{defn}


   The  reader unfamiliar with this orientation convention should verify the following observations: 

\begin{enumerate}
\item The above convention agrees with the usual convention for orienting products; (so given positive top-dimensional forms $\theta_{i}$ on $A_{i}$, $\theta_{1}\wedge\theta_{2}$ is a positive form on $A_{1}\times A_{2}$.)
\item Given two transverse submanifolds $A$ and $B$ of a manifold $M$, with normal bundles $N_{A}$ and $N_{B}$ oriented by the convention \[TA\oplus N_{A}=TM\text{ and }TB\oplus N_{B}=T{M}\] then $A\cap B$ considered as a fiber-product is oriented so that \[T (A\cap B)\oplus N_{B}\oplus N_{A}=TM\]  Be warned that some readers may consider this the usual convention for orienting $B\cap A$! See also example \ref{intersection}.
\item The orientation of $\ex B\fp gf\ex A$ differs from the orientation of $\ex A\fp fg\ex B$ by \[(-1)^{(\dim \ex A-\dim\ex C)(\dim\ex B-\dim\ex C)}\ .\]
To see this, note that swapping the order of the direct sum of $\coker Tf$  and $T\ex A$ and $\ker Tf$ and $TC$  gives \[\dim A\dim(\coker Tf)+\dim C\dim(\ker Tf)\] sign changes. Similarly, changing the orientation convention for the kernel relative to the cokernel of $Tg$ gives \[\dim B\dim(\coker Tg)+\dim C\dim(\ker Tg)\] sign changes. Then rearranging  $\coker Tf \oplus T\ex A\fp fg \ex B \oplus \coker Tg$ gives \[\begin{split}(\dim A+\dim B-\dim C)&(\dim(\coker Tf)+\dim(\coker Tg))\\&+\dim(\coker Tf)\dim(\coker Tg)\end{split}\] sign changes and rearranging   $\ker Tf\oplus T\ex C\oplus \ker Tg$ gives   \[\dim C(\dim(\ker Tf)+\dim(\ker Tg))+\dim (\ker Tf)\dim(\ker Tg)\] further sign changes. Summing these sign changes and simplifying gives the required expression.

\item The above convention for orienting the tangent space at a point of $\ex A\fp fg \ex B $ does give a well defined orientation on $\ex A\fp fg \ex B$. (You must check that deforming $Tf$ and $Tg$ continuously doesn't lead to any discontinuous change in orientation convention.)
\item If the normal bundle of $\ex A\fp fg\ex B\subset\ex A\times \ex B$ is identified with the pullback of $T\ex C$ using $Tf-Tg$, then 
the identification
\[T(\ex A\fp fg\ex B)\oplus T\ex C=T(\ex A\times \ex B)\]
changes orientation by the sign
\[(-1)^{\dim B\dim C}\] 
Of course, if we used $Tg-Tf$ to identify our normal bundle with the pull back of  $T\ex C$, then the sign would be $(-1)^{\dim C(\dim C+\dim B)}$, which agrees with the convention found on page 114 of \cite{hp}. 

\item The above convention makes the fiber-product associative in the sense that where defined, 
\[\ex A\fp f{g\circ k'} \lrb{\ex B\fp hk \ex C}=\lrb{\ex A\fp fg \ex B}\fp {h\circ f'}k\ex C=\lrb{\ex A\fp fg \ex B}\fp {f'}{k'}\lrb{\ex B\fp hk\ex C}\ .\]
The proof of associativity is not entirely trivial--- a sketch is below. It helps to consider the following commutative diagram:
{\setlength\mathsurround{0pt}
\[\begin{tikzcd}(\ex A\fp fg \ex B)\fp {f'}{k'}(\ex B\fp hk \ex C)\rar{f''}\dar{k''} &\ex B\fp hk\ex C
\rar{h'}\dar{k'}&\ex C\dar{k}
\\ \ex A\fp fg\ex B \rar{f'}\dar{g} &\ex B \rar{h}\dar{g} & \ex M_{2}
\\ \ex A\rar{f}&\ex M_{1}\end{tikzcd}\]}
Note that $Tf$, $Tf'$ and $Tf''$ have the same kernel and cokernel. Our orientation convention is equivalent to requiring that the relative orientations of these kernels and cokernels are the same, and that the orientation of $\ex A\fp fg \ex B$ is such that the following identifications are oriented isomorphisms.
\[\coker Tf'\oplus T\ex A\fp fg\ex B=\ker Tf'\oplus T\ex B\]
\[\coker Tf\oplus T\ex A=\ker Tf\oplus T\ex M_{1}\]
It follows that 
\[\ex A\fp f{g\circ k'} \lrb{\ex B\fp hk \ex C}=\lrb{\ex A\fp fg \ex B}\fp {f'}{k'}\lrb{\ex B\fp hk\ex C}\ .\]

Similarly, our orientation convention can be described only considering the downward pointing maps. All is as above, except that the kernels and cokernels now go on the right in the above identifications. It follows that
\[\lrb{\ex A\fp fg \ex B}\fp {h\circ f'}k\ex C=\lrb{\ex A\fp fg \ex B}\fp {f'}{k'}\lrb{\ex B\fp hk\ex C}\ .\]

\end{enumerate}

\

\begin{lemma} \label{fp integration} Suppose that $\ex A$, $\ex B$ and $\ex C$ are oriented exploded manifolds, \[f\co \ex A\longrightarrow \ex {C}\] is a $\Cs\infty1$ map, and  \[g\co \ex B\longrightarrow \ex {C}\] satisfies the conditions enumerated in Theorem \ref{fiber integration} for $g_{!}$ to exist. 

Consider the following commutative diagram involving the fiber-product of $f$ and $g$.
{\setlength\mathsurround{0pt}
\[\begin{tikzcd}\ex A\fp f{g}\ex {B}\rar{{g}'}\dar{f'}&\ex A\dar{f}
\\ \ex {B}\rar{{g}}&\ex {C}
\end{tikzcd}
\]}
Then ${g}'_{!}$ also exists, and the following diagram  is commutative. 
{\setlength\mathsurround{0pt}
\[\begin{tikzcd}\Omega_{c}^{*}(\ex A\fp f{g}\ex {B})\rar{{g}'_{!}}&\Omega_{c}^{*}(\ex A)
\\ (f'^{*})^{-1}\Omega_{c}^{*}(\ex A\fp f{g}\ex {B})\subset \Omega^{*}_{c}(\ex {B})\rar{{g}_{!}}\uar{f'^{*}}&\Omega_{c}^{*}(\ex {C})\uar{f^{*}}
\end{tikzcd}
\]}

\end{lemma}
\pf

As noted in the proof of Lemma 10.4 in \cite{iec}, if ${g}$ is a submersion and $T{g}$ is also surjective on integral vectors, ${g}'$ is a submersion which is also surjective on integral vectors. 
Therefore ${g}'$ satisfies the conditions of Theorem \ref{fiber integration} and ${g}'_{!}$ exists. 

It remains to verify that $f^{*}\circ {g}_{!}= {g}'_{!}\circ f'^{*}$ when restricted to forms $\theta\in\Omega^{*}_{c}(\ex B)$ so that $f'^{*}\theta$ has complete support in $\Omega^{*}(\ex A\fp fg\ex B)$. Note that $f'$ gives an isomorphism between the fibers of ${g}'$ and the fibers of ${g}$. When we consider $f'^{*}\theta$ as a top form $(f'^{*}\theta)'$ on the fibers with values in $g'^{*}\bigwedge T^{*}\ex A$, this form can be obtained from the corresponding form $\theta'$ by applying $f'^{*}\otimes ({g}'^{*}\circ f^{*})$. 
So
\[(f'^{*}\theta)'=(f'^{*}\otimes ({g}'^{*}\circ f^{*}))(\theta')\in \Omega^{\text{top}}_{c}({g}'^{-1}(p))\otimes {g}'^{*}\bigwedge T^{*}_{p}\ex A\ .\]

As we obtain ${g}'_{!}f'^{*}\theta$ by integrating $(g'^{*}\theta)'$ along the fibers of ${g}'$,  if the fibers of ${g}'$ are oriented the same as the fibers of ${g}$, 
\[f^{*}\circ {g}_{!}(\theta)={g}'_{!}\circ f'^{*}(\theta)\ .\]
Recall that to define integration along fibers, we orient so that 
\[T\ex {C}\oplus \ker T{g}=T\ex {B}\]
and 
\[T\ex A\oplus \ker T{g}'= T(\ex A\fp f{g}\ex {B})\]
 On the other hand, to orient $\ex A\fp f{g}\ex {B}$, we make the oriented identifications
 \[\coker Tf\oplus T\ex A=\ker Tf\oplus T\ex {C}\]
 \[T\ex {B}=T\ex {C}\oplus \ker T{g}\]
 \[\coker Tf\oplus T (\ex A\fp f{g}\ex {B})=\ker Tf\oplus T\ex {C}\oplus \ker T{g}\]
 Inserting the first of the above three equations into the last equation then gives that 
 \[T(\ex A\fp f{g}\ex {B})=T\ex A\oplus \ker T{g}\]
 is an oriented isomorphism.
 Therefore, with our orientation convention, the fibers of ${g}$ and ${g}'$ have the same orientation. It follows that
\[f^{*}\circ {g}_{!}(\theta)={g}'_{!}\circ f'^{*}(\theta)\]
as required.
\stop

  \section{Poincar\'e duality and fiber-products}\label{pdfp}
  
  In this section, we consider the relationship between Poincar\'e duality and fiber-products. In particular,  the relationship between Poincar\'e duality and both refinements  and intersection products.
  
  \
  
  Suppose that $\ex B'\longrightarrow \ex B$ is a refinement map (defined in section 10 of \cite{iec}). The corresponding map $H^{*}(\ex B)\longrightarrow H^{*}(\ex B')$ need not be an isomorphism. For example, suppose that $\ex B$ is a refinement of $\ex T^{n}$ corresponding to subdividing $\mathbb R^{n}$ into a the  toric fan of a nonsingular toric manifold. Then $H^{*}(\ex B)$ is isomorphic to the cohomology of the corresponding toric manifold $\totl{\ex B}$. Further subdividing this toric fan will produce a toric manifold with higher dimensional cohomology. 
  
  Suppose now that we have a map $f\co \ex C\longrightarrow \ex B$ where $\ex C$ is a complete oriented exploded manifold and Poincar\'e duality holds for $\ex B$. Then there exists some closed form $\theta\in \Omega^{*}_{c}(\ex B)$ so that for all $\alpha\in \Omega^{*}(\ex B)$, 
  \[\int_{\ex C}f^{*}\alpha=\int_{\ex B}\alpha\wedge \theta \] 
  This form $\theta$ may be unsatisfactory for the following reason: The fiber-product of $f$ with any refinement $\ex B'\longrightarrow \ex B$ gives a refined map $f'\co \ex C'\longrightarrow \ex B'$. Ideally, the pull back of $\theta$ to $\ex B'$ will then be the Poincar\'e dual to $f'$, but this may not be the case because there may be classes in $H^{*}(\ex B')$ which are not pulled back from classes in $H^{*}(\ex B)$. 
  
  \
  
 \begin{lemma}\label{poincare pullback} Let $\ex A$, $\ex B$ and $\ex C$ be oriented exploded manifolds. Suppose that $\ex C$ is  complete  and $f\co \ex C\longrightarrow \ex B$ is a $\Cs\infty1$ map so that  
 \[Tf\co {}^{\mathbb Z}T_{x}\ex C\longrightarrow {}^{\mathbb Z}T_{f(x)}\ex B\]
 is surjective. Then given any neighborhood $N$ of $f(\ex C)\subset \ex B$, there exists a closed form $\theta\in\Omega^{*}_{c}(\ex B)$, supported in $N$, and  Poincar\'e dual to $f$ in the sense that 
 \[\int_{\ex B}\alpha\wedge \theta=\int_{\ex C}f^{*}\alpha \ \ \ \ \ \ \text{ for all closed }\alpha\in \Omega^{*}(\ex B).\]
 Suppose that $g\co \ex A\longrightarrow \ex B$ is any complete $\Cs\infty1$ map transverse to $f$. Then $g^{*}\theta$ is Poincar\'e dual to the map $f'$ below  
{\setlength\mathsurround{0pt} \[\begin{tikzcd}\ex A\fp gf\ex C\rar{g'}\dar{f'} &\ex C\dar{f}
 \\ \ex A\rar{g}&\ex B\end{tikzcd}\]}
 
 in the sense that 
 \[\int_{\ex A}\alpha\wedge g^{*}\theta=\int_{\ex A\fp gf \ex C}f'^{*}\alpha
 \ \ \ \ \text{ for all closed }\alpha\in \Omega^{*}(\ex A).\]
 
 \end{lemma} 
 
\pf 

  Extend $f$ to a submersion $h\co \ex C\times \mathbb R^{n}\longrightarrow \ex B$ satisfying the conditions of Theorem \ref{fiber integration}. (Here $h$ extends $f$ in the sense that  $h(p,0)=f(p)$.) Choose a compactly supported form $\theta_{0}$ on $\mathbb R^{n}$ that integrates to $1$, consider this form $\theta_{0}$ as a form on $\ex C\times \mathbb R^{n}$, then integrate along the fibers of $h$ to obtain 
  \[\theta:=h_{!}\theta_{0}\in \Omega^{*}_{c}(\ex B)\]
  
   This form $\theta$ represents the Poincar\'e dual of $f$. In particular, suppose that $\alpha\in \Omega^{*}(\ex B)$ is closed. Then our adaptation of Stokes' theorem,  Theorem \ref{stokes}, implies that
  \[\int_{\ex C\times x\subset \ex C\times \mathbb R^{n}}h^{*}\alpha=\int_{\ex C\times 0}h^{*}\alpha=\int_{\ex C}f^{*}\alpha\ .\] 
  Therefore, 
  \[\int_{\ex B}\alpha\wedge h_{!}\theta_{0}=\int_{\ex C\times \mathbb R^{n}}h^{*}\alpha\wedge\theta_{0}=\int_{\ex C}f^{*}\alpha\ .\]
  So $\theta=h_{!}\theta_{0}$ is Poincar\'e dual to $f$. By choosing $\theta_{0}$ supported close to $0\in\mathbb R^{n}$ we may arrange that $\theta$ is supported close to the image of $f$.
  
 Given our complete map $g$ transverse to $f$, we may now consider the following fiber-product.
 {\setlength\mathsurround{0pt}
 \[\begin{tikzcd}\ex A\fp gh(\ex C\times \mathbb R^{n})\rar{g'}\dar{h'} &\ex C\times \mathbb R^{n}\dar{h}
 \\ \ex A\rar{g}&\ex B\end{tikzcd}\]}
 
 Applying Lemma \ref{fp integration} gives 
 \[g^{*}\theta=g^{*}(h_{!}(\theta_{0}))=h'_{!}(g'^{*}(\theta_{0}))\]
 so 
 \[\int_{\ex A}\alpha\wedge g^{*}\theta=\int_{A\fp gh(\ex C\times \mathbb R^{n})}h'^{*}\alpha\wedge g'^{*}(\theta_{0})\ .\]
 
 Define the map $r\co [0,1]\times\ex C\times \mathbb R^{n}\longrightarrow \ex B$ by
 \[r(t,p,x)=h(p,xt)\]
As $f(p)=h(p,0)$ and  $f$ and $h$ are transverse to $g$, our new map $r$ is also transverse to $g$, so we may take the following fiber-product.
 {\setlength\mathsurround{0pt} \[\begin{tikzcd}\ex A\fp gr([0,1]\times\ex C\times \mathbb R^{n})\rar{\hat g'}\dar{r'} & {[0,1]}\times\ex C\times \mathbb R^{n}\dar{r}
 \\ \ex A\rar{g}&\ex B\end{tikzcd}\]}
   As $g$ is complete, $\hat g'$ is also complete, so $\hat g'^{*}\theta_{0}$ is completely supported. We may now apply Stokes' theorem. Our map $r$ restricted to $t=1$ is $h$, and restricted to $t=0$ is 
  \[r(0,p,x)=f(p)\ .\] 
  Associativity of fiber-products implies that the corresponding boundary of $\ex A\fp gr([0,1]\times\ex C\times \mathbb R^{n})$ is equal to $(A\fp gf\ex C)\times \mathbb R^{n}$. Then

  \[\int_{A\fp gh(\ex C\times \mathbb R^{n})}h'^{*}\alpha\wedge g'^{*}\theta_{0}=\int_{(A\fp gf\ex C)\times \mathbb R^{n}} f'^{*}\alpha\wedge \theta_{0}\]
  where in the above, $f'^{*}\alpha$ and $\theta_{0}$ indicate the pullback of the corresponding forms on $\ex A\fp gf\ex C$ and  $\mathbb R^{n}$ respectively. Therefore, 
  \[\int_{A}\alpha\wedge g^{*}\theta=  \int_{(A\fp gf\ex C)\times \mathbb R^{n}}f'^{*}\alpha\wedge \theta_{0}=\int_{A\fp gf\ex C}f'^{*}\alpha\]
as required.

\stop

\

\begin{example}[Intersection of submanifolds and Poincar\'e duality]\label{intersection}\end{example}
 Suppose that $\ex A$ and $\ex C$ are complete  exploded manifolds and oriented submanifolds of the oriented exploded manifold $\ex B$ in the sense that they can be locally described as the inverse image of a regular value of some $\mathbb R^{n}$--valued $\Cs\infty1$  function. Then we may use the construction of Lemma \ref{poincare pullback} to construct Poincar\'e duals $\theta_{\ex A}$ and $\theta_{\ex C}$ to $\ex A$ and $\ex C$. If $\ex A$ and $\ex C$ are transverse, then Lemma \ref{poincare pullback} implies that the $\theta_{\ex C}$ restricted to $\ex A$ is Poincar\'e dual to $\ex A\cap \ex C\subset \ex A$.

Therefore, 
\[\int_{\ex A\cap \ex C}\alpha=\int_{\ex A}\alpha\wedge\theta_{\ex C}=\int_{\ex B}\alpha\wedge \theta_{\ex C}\wedge\theta_{\ex A}\]
and the Poincar\'e dual to $\ex A\cap \ex C$ is $\theta_{\ex C}\wedge \theta_{\ex A}$. So, with our sign convention, intersection products correspond under Poincar\'e duality to wedge products with the order reversed.

\

Be warned that if  $\ex A$ and $\ex C$ are not submanifolds in the above sense, the above formula may not hold. For example, let $\ex B$ be a refinement of $\ex T^{2}$ corresponding to dividing $\mathbb R^{2}$ into the standard quadrants, and consider $\ex A:=\{\tilde z_{1}=\tilde z_{2}\}\subset \ex B$ and $\ex C\co \{(\tilde z_{1}+\tilde z_{2}+1\e1)\in 0\e{\mathbb R}\}\subset \ex B$. Note that for any $\theta\in\Omega^{*}(\ex B)$
\[\int_{\ex C}\theta=\int_{\ex C'}\theta\]
where $\ex C':=\{\tilde z_{1}=-\tilde z_{2}\}$. This is because $\theta$ must vanish out where $\ex C$ and $\ex C'$ differ.
 Therefore, the Poincar\'e duals of $\ex C$ and $\ex C'$ are the same, so if the usual relationship between intersections and wedge products held, the Poincar\'e dual of  $\ex A\cap \ex C$ should be equal to the Poincar\'e dual of $\ex A\cap \ex C'$. But $\ex A\cap \ex C$ is a single point and $\ex A\cap \ex C'$ is empty.  
 
 The solution to this problem is to allow a more flexible class of differential forms called refined forms.

 \section{Refined cohomology}

 \begin{defn} A refined form $\theta\in\ro^{*}(\ex B)$ is a choice of $\theta_{p}\in\bigwedge T^{*}_{p}(\ex B)$ for all $p\in\ex B$ so that given any point $p\in \ex B$, there exists an open neighborhood $U$ of $p$ and a complete, surjective, equidimensional submersion
 \[r\co U'\longrightarrow U\]
   so that there is a form $\theta'\in\Omega^{*}(U')$ pulling back  $\theta$ in the sense  that if $v$ is any vector on $U'$ so that $dr(v)$ is a vector based at $p$, then 
  \[\theta'(v)=\theta_{p}(dr(v))\ .\] 
 
 A refined form  $\theta\in\ro^{*}(\ex B)$ is completely supported if there exists some complete subset $V$ of an exploded manifold $\ex C$ with a map $\ex C\longrightarrow \ex B$ so that $\theta_{p}=0$ for all $p$ outside the image of $V$. Use the notation $\ro^{*}_{c}$ for completely supported refined forms. 
 
 Denote the homology of $(\ro^{*}(\ex B),d)$ by $\rh^{*}(\ex B)$ and  $(\ro^{*}_{c}(\ex B),d)$ by $\rh^{*}_{c}(\ex B)$.
 \end{defn}
 
In Theorem \ref{rfi} we shall show that refined forms push forward along  oriented submersions, removing the condition on integral vectors from Theorem \ref{fiber integration}. This allows us to 
show in Lemma \ref{rpp} that the Poincar\'e dual to any oriented map from a complete exploded manifold exists as a refined differential form, and that the familiar relationship between fiber-products of maps and wedge products of these Poincar\'e duals holds. In fact,  refined cohomology is minimal extension of $H^{*}$ with pushforwards and  Poincar\'e duality compatible with fiber products. 

\

For defining the refined cohomology above, it should be obvious that \[d\co \ro^{k}(\ex B)\longrightarrow \ro^{k+1}(\ex B)\]is well defined and $d^{2}=0$. Less immediate, is the fact that $\ro^{*}(\ex B)$ is closed under addition and wedge products.  If $\theta_{1}$ and $\theta_{2}$ are refined forms, then any point  $p$ has a neighborhood $U$ with complete, surjective, equidimensional submersions
\[r_{i}\co U'_{i}\longrightarrow U\]
so that $r_{i}^{*}\theta_{i}\in\Omega^{*}(U'_{i})$. Taking the fiber-product of $r_{1}$ with $r_{2}$ gives a complete, surjective, equidimensional submersion
\[r'\co U_{1}\fp{r_{1}}{r_{2}}U_{2}\longrightarrow U\]
 so that $r'^{*}\theta_{i}\in \Omega^{*}(U_{1}\fp{r_{1}}{r_{2}}U_{2})$. Therefore, $\theta_{1}+\theta_{2}$ and $\theta_{1}\wedge\theta_{2}$ are  in $\ro^{*}(\ex B)$. 
 
 \

 The existence of partitions of unity combined with Lemma \ref{finite integral} implies that the integral of $\theta\in\ro_{c}^{*}(\ex B)$ over $\ex B$ is finite and well defined. In particular if $\rho\theta$ is supported in $U$ and the map $r\co U'\longrightarrow U$ has degree $m$, then 
 \[\int_{U}\rho\theta:=\frac 1m\int_{U'}r^{*}\rho\theta\ .\]

\

Note also that given any $\Cs\infty1$ map $f\co \ex A\longrightarrow \ex B$, there is a linear chain map 
\[f^{*}\co \ro^{*}(\ex B)\longrightarrow \ro^{*}(\ex A)\]
 defined as usual so that 
\[(f^{*}\theta)_{p}(v):=\theta_{f(p)}(Tf(v))\ .\]
To see that $f^{*}\theta$ is actually in $\ro^{*}(\ex A)$, let $r\co U'\longrightarrow U$ be a complete, 
equidimensional submersion onto a neighborhood of $f(p)$ so that $r^{*}\theta\in \Omega^{*}(U')$. Then   
taking the fiber-product of $r\co U'\longrightarrow \ex B$ with $f$ gives a complete equidimensional submersion onto a neighborhood of $p$ so that the pullback of $f^{*}\theta$ is in $\Omega^{*}$, so $f^{*}\theta\in \Omega^{*}(\ex A)$. 
 
 \
 
  Our version of Stokes' theorem also extends  to refined forms in $\ro^{*}_{c}(\ex B)$. If $\ex B$ is a complete exploded manifold, the integration pairing on $\rh^{*}(\ex B)$ is nondegenerate, but as $\rh^{*}(\ex B)$ is in general infinite-dimensional, this does not imply Poincar\'e duality.
 
 \

 \begin{thm}\label{rfi} Given any submersion $f\co \ex B\longrightarrow\ex C$ between oriented exploded manifolds, there exists a linear chain map 
 \[f_{!}\co \ro^{*}_{c}(\ex B)\longrightarrow \ro^{*}_{c}(\ex C)\]
 uniquely determined by the property that
 \[\int_{\ex C}\alpha\wedge  f_{!}\beta=\int_{\ex B}(f^{*}\alpha)\wedge\beta\]
for all $\beta\in\ro^{*}_{c}(\ex B)$ and $\alpha\in\ro^{*}(\ex C)$.
 \end{thm}
 
 \pf
 
Given any point $p\in \ex C$, we may take a refinement of a neighborhood of $p$ so that the inverse image of $p$ in the refined neighborhood is contained in a stratum which is a smooth manifold. As a smooth form on a manifold is determined by its integral against compactly supported forms, $f_{!}\beta$ around $p$ is uniquely determined by the property $\int_{\ex C}\alpha\wedge  f_{!}\beta=\int_{\ex B}(f^{*}\alpha)\wedge\beta$. As the right hand side of this equation is linear in $\beta$, it follows that $f_{!}$ is linear if it exists. Stokes' theorem implies that if $\alpha\in \ro^{k}(\ex C)$, 
\[\int_{\ex C}\alpha\wedge df_{!}\beta=(-1)^{k+1}\int_{\ex C}(d\alpha)\wedge f_{!}\beta=(-1)^{k+1}\int_{\ex B}(df^{*}\alpha)\wedge\beta=\int_{\ex B}(f^{*}\alpha)\wedge d\beta\] 
  so if $f_{!}$ exists, it is a linear chain map. Using a partition of unity, we may restrict to the case that $f$ is a map between coordinate charts $U$ and $V$ and $\beta$ pulls back to a form in $\Omega^{*}_{c}(U')$. Then, using a partition of unity on $U'$, we may restrict to the case that $\beta$ is supported in a single coordinate chart of $U'$. By relabeling we do not lose generality by assuming that $\beta$ is supported in a single coordinate chart $U$. 

The tropical part of $U$ and  $V$ are polytopes $\totb U$ and $\totb{V}$. There exists a coordinate chart $V'$ with a complete equidimensional submersion $V'\longrightarrow V$ so that the image of integral vectors from $U$ in $V$ is always a full sublattice of the image of integral vectors from $V'$. (The tropical part of $V'$ is $\totb V$ with a different integral-affine structure.) Then we may choose a refinement $V''\longrightarrow V$ corresponding to a subdivision of $V$ so that $\totb{f(U)}$ is a polytope in this subdivision. Suppose that $\alpha\in \ro^{*}V$ pulls back to a $\Cs\infty1$ form on some $V'''\longrightarrow V$. Then let $\hat V$ be the fiber-product over $V$ of $V'$ with $V''$ and $V'''$,
and let $r\co \hat U\longrightarrow U$ be the fiber-product of $\hat V\longrightarrow V$ with $f\co U\longrightarrow V$.   

 Then $\hat f\co \hat U\longrightarrow \hat V$ is a submersion which also is surjective on integral vectors. Therefore,  Theorem \ref{fiber integration} implies that there is a linear chain map $\hat f_{!}\co \Omega^{*}_{c}(\hat U)\longrightarrow \Omega^{*}_{c}(\hat V)$ so that 
 \[\int_{\hat V}\alpha\wedge \hat f_{!}(r^{*}\beta)=\int_{\hat U}(\hat f^{*}\alpha)\wedge r^{*}\beta\]
 for all $\alpha\in \Omega^{*}(\hat V)$ and $\beta\in\Omega^{*}_{c}(U)$. Considering $f_{!}\beta$ as refined form in $\ro^{*}_{c}(\ex B)$, we have our map $f_{!}$. As $\hat U\longrightarrow U$ has the same degree as $\hat V\longrightarrow V$, the above formula implies that
 \[\int_{V}\alpha\wedge f_{!}\beta=\int_{U}(f^{*}\alpha)\wedge \beta\]

   Lemma \ref{fp integration} implies that this map $f_{!}$ is independent of further refinement of $\hat V$ and $\hat U$, so $f_{!}\beta$ depends only on $\beta$ as an element of $\ro^{*}_{c}(U)$, not on $\beta$ as an element of $\Omega^{*}_{c}(\hat U)$. 
 
    \stop

 \

 \begin{lemma} \label{rfpi} Suppose that $\ex A$, $\ex B$ and $\ex C$ are oriented exploded manifolds, $f\co \ex A\longrightarrow \ex B$ is a submersion, and $g\co \ex C\longrightarrow \ex B$ is a $\Cs\infty1$ map. 

Consider the following commutative diagram involving the fiber-product of $g$ and $f$.
{\setlength\mathsurround{0pt}\[\begin{tikzcd}\ex C\fp gf\ex A\rar{f'}\dar{g'}&\ex C\dar{g}
\\ \ex A\rar{f}&\ex B
\end{tikzcd}
\]}
Then $f'_{!}$ also exists, and the following diagram  is commutative. 
{\setlength\mathsurround{0pt}
\[\begin{tikzcd}\ro_{c}^{*}(\ex C\fp gf\ex A)\rar{f'_{!}}&\ro_{c}^{*}(\ex C)
\\ (g'^{*})^{-1}\lrb{\ro_{c}^{*}(\ex C\fp gf\ex A)}\subset\ro^{*}_{c}(\ex A)\uar{g'^{*}}\rar{f_{!}}&\ro_{c}^{*}(\ex B)\uar{g^{*}}
\end{tikzcd}
\]}

\end{lemma}

\pf This lemma has the same proof as Lemma \ref{fp integration}, except Theorem \ref{rfi} is used instead of Theorem \ref{fiber integration}.

\stop
 
 \
 
 \begin{remark} \label{relative orientation} For any submersion, $f\co \ex A\longrightarrow \ex B$, between oriented exploded manifolds, giving $\ex A$ and $\ex B$ the opposite orientations results in the same $f_{!}$.  (One way to see this from Lemma \ref{rfpi} is to set $\ex C$  to be $\ex B$ with the opposite orientation.) It follows that $f_{!}$ only depends on the choice of relative orientation of $f$. (For a given map $f\co \ex A\longrightarrow \ex B$, a relative orientation is a homotopy class of nonvanishing section of $\wedge^{top}(T\ex A\otimes f^{*}T\ex B)$.)

 In the case of a relatively oriented submersion $f\co \ex A\longrightarrow \ex C$ where $\ex C$ is not necessarily oriented, Lemma \ref{rfpi} implies that we may  define $f_{!}\co \ro_{c}^{*}\ex A\longrightarrow \ro_{c}^{*}\ex C$ to be the unique linear chain map satisfying the property that given any map $g\co \ex B\longrightarrow \ex C$ from an oriented exploded manifold $\ex B$, the diagram
 {\setlength\mathsurround{0pt} \[\begin{tikzcd}\ro_{c}^{*}(\ex C\fp gf\ex A)\rar{f'_{!}}&\ro_{c}^{*}(\ex C)
\\ (g'^{*})^{-1}\lrb{\ro_{c}^{*}(\ex C\fp gf\ex A)}\subset\ro^{*}_{c}(\ex A)\uar{g'^{*}}\rar{f_{!}}&\ro_{c}^{*}(\ex B)\uar{g^{*}}
\end{tikzcd}
\]}
commutes, where $f'_{!}$ is defined by Theorem \ref{rfi} because the relative orientation of $f$ induces a relative orientation of $f'$, which together with the orientation of $\ex B$ gives an orientation for $(\ex C\fp gf\ex A)$. 
 \end{remark}
 
 Even though the perfect pairing version of Poincar\'e duality does not necessarily hold for $\rh^{*}(\ex B)$, the following lemma gives an analogue of the Poincar\'e dual of a map from a complete manifold. 
 
 \begin{lemma}\label{rpp} Suppose that $\ex C$ is a complete oriented exploded manifold and $f\co \ex C\longrightarrow \ex B$ is a $\Cs\infty1$ map to an oriented exploded manifold $
\ex B$. Then given any metric on $\ex B$ and distance $r$,  there exists a closed form $\theta\in\ro^{*}_{c}(\ex B)$ supported within a radius $r$ of $f(\ex C)$ which is Poincar\'e dual to $f$ in the sense that 
 \[\int_{\ex B}\alpha\wedge \theta=\int_{\ex C}f^{*}\alpha \ \ \ \ \ \ \text{ for all closed }\alpha\in \ro^{*}(\ex B).\]
 Suppose that $\ex A$ is oriented and $g\co \ex A\longrightarrow \ex B$ is any complete $\Cs\infty1$ map transverse to $f$. Then $g^{*}\theta$ is Poincar\'e dual to the map $f'$ below  
 {\setlength\mathsurround{0pt}\[\begin{tikzcd}\ex A\fp gf\ex C\dar{f'}\rar{g'} &\ex C\dar{f}
 \\ \ex A\rar{g}&\ex B\end{tikzcd}\]}
 in the sense that 
 \[\int_{\ex A}\alpha\wedge g^{*}\theta=\int_{\ex A\fp gf \ex C}f'^{*}\alpha
 \ \ \ \ \text{ for all closed }\alpha\in \ro^{*}(\ex A).\]
 
 \end{lemma} 
 
 \pf 
 
 The proof of this lemma is identical to the proof of Lemma \ref{poincare pullback}, except Theorem \ref{rfi} is used in the place of Theorem \ref{fiber integration}, and Lemma \ref{rfpi} is used instead of Lemma \ref{fp integration}.
 
 \stop
 
 \section{Partitions of unity and good covers}\label{partition}

Throughout this paper, we are assuming that our exploded manifolds considered as topological spaces are second countable. The following lemma constructs a partition of unity subordinate to a given open cover of an exploded manifold.

\begin{lemma} Given any open cover $\{U_{\alpha}\}$ of an exploded manifold $\ex B$, there exists a partition of unity subordinate to $\{U_{\alpha}\}$.
\end{lemma}

\pf

Any (second countable) exploded manifold has an exhaustion by compact subsets $K_{i}$ so that $K_{i-1}$ is contained in the interior of $K_{i}$. This fact follows from the observation that it holds for  $\mathbb R^{n}\times \et mP$, and any (second countable) exploded manifold has a countable cover by open subsets isomorphic to $\mathbb R^{n}\times \et mP$. 

A second ingredient needed for construction of partitions of unity is the existence of bump functions. There exists a smooth function with compact support which is positive on any given compact subset of $\mathbb R^{n}\times \et mP$. Given any point  $p$ in an open subset $U$ of an exploded manifold, Lemma 6.10 of \cite{iec} states that there exists an open neighborhood of $p$ contained inside $U$ which is isomorphic to $\mathbb R^{n}\times\et mP$. Therefore, there exists a smooth nonnegative function which is positive at $p$ and which has support compactly contained inside $U$.

We may now construct partitions of unity as usual. Let  $\{U_{\alpha}\}$ be any open cover  of $\ex B$.  For each point $p$ in $K_{i}\setminus K_{i-1}$, choose a nonnegative bump function $\rho_{p}$ which is positive at $p$ and which has compact support contained inside $K_{i+1}\setminus K_{i-1}$ and some $U_{\alpha}$.  The sets $\{\rho_{p}>0\}$ form an open cover of $\ex B$ which have a locally finite subcover $\{\rho_{p_{i}}>0\}$ for $i=1,\dotsc$.  Then $\sum_{i} \rho_{p_{i}}$ is smooth and positive, so we may divide our functions  $\rho_{p_{i}}$ by this sum to obtain the required partition of unity.  

\stop

\

\begin{lemma} \label{good cover}Any compact exploded manifold $\ex B$ has a finite good cover $\{U_{i}\}$ in the sense that the intersection of any number of these $U_{i}$ is either empty or isomorphic to $\mathbb R^{n}\times\et mP$.

\end{lemma}

\pf 

One way to prove that a manifold $M$ has good cover is to choose a connection and then construct a cover of $M$ by geodesically convex open sets. This proof does not immediately generalize to exploded manifolds,  because geodesically convex open neighborhoods of points in strata with nonzero tropical dimension do not exist. In this proof, we shall first choose a nice `equivariant' set of coordinate charts, construct a connection $\nabla$ compatible with these coordinate charts, construct functions which are convex with respect to $\nabla$, then use these convex functions to construct an open cover satisfying a convexity condition strong enough to prove that it is a good cover.

\

It was shown in \cite{cem} that any exploded manifold has a cover by equivariant coordinate charts isomorphic to open subsets of $\mathbb R^{n}\times \et mP$. (This was proved for smooth exploded manifolds in Lemma A3 of \cite{cem}. The same proof works for $\Cs\infty1$ exploded manifolds: the only modification necessary is that $\Cs\infty1$ coordinate charts and vectorfields should be used in place of smooth coordinate charts and vectorfields in the proof of Lemma A3.) Using equivariant coordinate charts means that each transition map or its inverse is in the form of a map

\begin{equation}\label{equivariant def}(x,\tilde z)\mapsto (f(x),g_{1}(x)\tilde z^{\alpha^{1}},\dotsc,g_{n}(x)\tilde z^{\alpha^{n}})\end{equation}

In particular, transition maps of the above type send the lattice of vectorfields $N$ generated by the real and imaginary parts of $\tilde z_{i}\frac\partial{\partial \tilde z_{i}}$ to a sublattice of the corresponding lattice in the target. Note that the coordinate charts with more structure are those with higher dimensional tropical part, so these equivariant transition maps never decrease the dimension of the tropical part. We shall  assume that we have a finite number of coordinate charts and that if two coordinate charts intersect, the tropical part of one of the coordinate charts is a face of the other coordinate chart. Recall that the closure of any stratum is a face - the above statement does not assume that it is a codimension 1 face!)

\

By using a partition of unity and reducing the size of our equivariant coordinate charts where necessary, we can choose a connection $\nabla$ on $T\ex B$ so that in our coordinate charts, for any vector $w$ in the lattice of vectorfields $N$ generated by the real and imaginary parts of $\tilde z_{i}\frac\partial{\partial \tilde z_{i}}$,
\[ \nabla_{w}=L_{w}\text{ and }\nabla w=0\]
To achieve the above, proceed as follows: in each coordinate chart, the standard flat connection obeys the above conditions. Choose a finite partition of unity consisting of bump functions compactly supported inside our equivariant coordinate charts, and average these flat connections using this partition of unity. Now reduce the size of coordinate charts so that the above condition holds, starting with the coordinate charts of the highest tropical dimension. Because our coordinate changes are equivariant, on the open subset of any of our  coordinate charts which is the complement of the support of all bump functions supported inside charts with lower tropical dimension, the averaged connection obeys the above conditions. We can safely reduce the size of our coordinate chart to this open set  because the complement is covered by charts with lower dimensional tropical part.

Our connection $\nabla$ on (an open subset of ) $\mathbb R^{n}\times\et mP$ defines a connection $\nabla'$ on an open subset of $\mathbb R^{n}$  follows: Let $x$ denote the projection to $\mathbb R^{n}$, and let $\tilde v$ indicate any lift of a vectorfield $v$ from $\mathbb R^{n}$ to $\mathbb R^{n}\times \et m{P^{\circ}}$ so that $Tx(\tilde v)=v$. Then  define \[\nabla'_{v_{1}}v_{2}:=Tx(\nabla_{\tilde v_{1}}\tilde v_{2})\] Note that every smooth vectorfield on $\mathbb R^{n}\times \et m{P^{\circ}}$ is independent of $\et m{P^{\circ}}$ as $P^{\circ}$ is an open polytope, so $Tx(\nabla_{\tilde v_{1}}\tilde v_{2})$ is indeed a vectorfield on $\mathbb R^{n}$. We'll check below that our conditions on $\nabla$ ensure that this projected connection is well defined: 

 Let $x_{i}$ be $x$ followed by projection of $\mathbb R^{n}$ onto the $i$th coordinate. If $w$ is in the kernel of $Tx$, then $\nabla_{w}(dx_{i}(\tilde v))=0$ if $\tilde v$ is lifted. As such a $w$ must be in the span of $N$, $\nabla_{w}dx_{i}=0$, so  $\nabla_{w}\tilde v$ is also in the kernel of $Tx$. This implies that $Tx(\nabla_{\tilde v_{1}}\tilde v_{2})$ does not depend on the choice of lift of $v_{1}$. The fact that $\nabla_{v}w=0$ for any $w\in N$ implies that if we instead choose $w$ in the kernel of $Tx$ (so $w$ is a sum of smooth functions times vectorfields in $N$),  then $\nabla_{v}w$ will also be in the kernel of $Tx$. It follows that $Tx(\nabla_{\tilde v_{1}}\tilde v_{2})$ does not depend on the choice of lift of $v_{2}$, and the connection $\nabla'$ is well defined.

\

We shall now locally construct some convex functions. In one of our equivariant coordinate charts, consider the function $\abs x^{2}$. When it is small enough, this function is (nonstrictly) convex in the sense that restricted to any $\nabla$--geodesic, it has nonnegative second derivative. In particular, if a geodesic has velocity $v$, then the second derivative of $\abs x^{2}$ restricted to the geodesic is 
\[\nabla_{v}(d\abs x^{2})(v)=2\sum_{i}(dx_{i}(v))^{2}+2x_{i}(\nabla_{v}dx_{i})(v)\ .\]
Restricted to vectors in the subspace generated by $\frac \partial {\partial{x_{i}}}$, the above  quadratic form is positive definite at $x=0$, and therefore positive definite on this subspace for $\abs x$ small enough. Let $w$ be any vector field given by a sum of constants times the real and imaginary parts of $\tilde z_{i}\frac\partial{\partial \tilde z_{i}}$, then 
\[\nabla_{v+w}(d\abs x^{2})(w)=\nabla_{v+w}(d\abs x^{2}(w))=0\ .\] 
Also, the fact that $\nabla_{w}=L_{w}$ implies that
\[\nabla_{w}(d\abs x^{2})=0\ .\] 
Therefore, 
\[\nabla_{v+w}(d\abs x^{2})(v+w)=\nabla_{v}(d\abs x^{2})(v)\ .\] 
As our quadratic form is positive definite on one subspace and independent of  a complimentary subspace, it is positive semidefinite.

Now  construct a proper convex function on (a subset of) our coordinate chart $\mathbb R^{n}\times \et mP$. Choose some basis $\{\zeta_{\alpha}:=\totl{\e a\tilde z^{\alpha}}\}$ of smooth monomials on $\et mP$, and consider the function 
\begin{equation}\label{fdef}f:=\sum_{\alpha}\abs{\zeta_{\alpha}^{2}}\ .\end{equation}
Restricted to geodesics in the directions spanned by $N$, this function is  (non strictly) convex, as verified by the following calculation: If $w_{i}$ indicates the real part of $\tilde z_{i}^{-1}d\tilde z_{i}(w)$, then
\[\nabla_{w}(df)(w)=\sum_{\alpha,i}4\abs{\zeta_{\alpha}}^{2}\alpha_{i}^{2}w_{i}^{2}\ .\] 
For $x$ small,  $\abs x^{2}$ is strictly convex  on the complementary subspace to $N$ spanned by $\frac\partial {\partial x_{i}}$. Therefore, if we choose $\lambda$ large enough,  $f+\lambda\abs x^{2}$ will be convex when it is $\leq 1$.

\

\begin{claim}\label{con cover}There exists a finite cover of $\ex B$ by open
 subsets $U_{i}$ of our coordinate charts given by \[U_{i}:=\left\{g_{i}<1\right\}\subset\mathbb R^{n}\times \et mP\] where 
 \[g_{i}\co \mathbb R^{n}\times \et mP\longrightarrow [0,\infty)\]
  is proper and (non-strictly) convex on the subset 
 \[2U_{i}:=\{g_{i}<2\}\subset \mathbb R^{n}\times \et mP\]
 so that every intersection $U$ of a finite number of these coordinate charts  satisfies the following two convexity assumptions:
 \begin{enumerate}
\item\label{con1} The projection $U'$ of $U\cap\mathbb R^{n}\times \et m{P^{\circ}}$ to $\mathbb R^{n}$ is geodesically convex using the connection $\nabla'$. 
\item\label{con2} $U$  is defined by some finite number of inequalities $g_{i}<1$ where each $g_{i}$ is a finite sum of positive functions on $U'$ times the square absolute value of smooth monomial functions, and $g_{i}$ is proper restricted to each $\et mP$ fiber of $U'\times \et mP$.  
\end{enumerate}
\end{claim}

We shall prove, in Claim \ref{con good} below, that the  above two convexity conditions imply that $U$ is isomorphic to $\mathbb R^{n}\times \et mP$. Condition \ref{con1} will allow us to choose an isomorphism of $U'$ with $\mathbb R^{n}$, and condition \ref{con2} will allow us to stretch $U$ in the $\et mP$ direction to give an isomorphism with $U'\times \et mP$. One important aspect of condition \ref{con2} is encapsulated by the following observation:

\begin{remark}\label{pullback monomial}Let $\psi\co U_{1}\longrightarrow U_{2}$ be an equivariant map between coordinate charts so that $\totb\psi$ is an inclusion of $\totb{U_{1}}$ as a face of $\totb{U_{2}}$. If   $\zeta$ is a smooth monomial function on $U_{2}$, 
\[\abs\zeta^{2}\circ \psi\]  
may be written as the square absolute value of some monomial function on $U_{1}$ times (the pullback to $U_{1}$ of) a smooth positive function on $U_{1}'$,
and if $x$ is (the pullback to $U_{2}$ of) a smooth positive function on $U_{2}'$, then 
\[x\circ \psi\] 
is also the pullback to $U_{1}$ of a smooth positive function on $U_{1}'$.
\end{remark}

Remark \ref{pullback monomial} follows directly from the the fact that our equivariant map $\psi$ may be written in the form of equation (\ref{equivariant def}).

\

 Now we choose a function $g$ satisfying the conditions of the $g_{i}$ from Claim \ref{con cover}. Setting $g$ to be a large enough multiple of $f+\lambda\abs{x-p}^{2}$ gives an open subset $U:=\{g<1\}$ satisfying the above convexity conditions: Condition \ref{con2} is obviously satisfied, and condition \ref{con1} is satisfied for $g$ a large enough multiple of $f+\lambda\abs{x-p}^{2}$ because  $f+\lambda\abs {x-p}^{2}$ restricted to $U'$ is the convex function $\lambda\abs{x-p}^{2}$. Within any coordinate chart isomorphic to an open subset of $\mathbb R^{n}\times \et mP$, we may therefore cover an open neighborhood of this chart intersected with $\mathbb R^{n}\times \et m{P^{\circ}}$ with open subsets individually satisfying the convexity conditions of Claim \ref{con cover}.

\

Suppose that we have a finite cover  of all strata of tropical dimension greater than $k$ by sets $U_{i}$ defined by functions $g_{i}$ satisfying  the above convexity conditions of Claim \ref{con cover}. We shall extend this cover to a cover of all strata of tropical dimension greater than or equal to  $k$ while still satisfying the convexity conditions of Claim \ref{con cover}. This will prove Claim \ref{con cover}.

 Choose a cover of the strata of tropical dimension $k$ using coordinate charts with tropical dimension $k$ so that each of these coordinate charts includes in an old coordinate chart via an equivariant map.  We may now cover the strata of dimension $k$ by open sets $U$ coming from functions $g$ defined in these new coordinate charts and satisfying the above convexity conditions, and satisfying the extra condition that if $U$ intersects a member $U_{i}$ of our previously constructed finite good cover, then $U$ is contained entirely inside $2U_{i}$. This new collection of open sets together with our old  cover is an open cover of the set of strata of dimension at least $k$, which is compact, so we can choose a finite sub cover. It remains to prove that this subcover satisfies the convexity conditions of Claim \ref{con cover}.  

The intersection of a finite number of these new sets $U$ satisfying the  two convexity conditions  of Claim \ref{con cover}  clearly still satisfies these convexity conditions, because all transition maps \emph{and} their inverses are equivariant. Intersection with some of the previously constructed $U_{i}$ then corresponds to restricting  to a subset where the functions $g_{i}<1$. Restricting a geodesically convex set to a set where a (nonstrictly) convex function  is less than $1$ gives a geodesically convex set, so this intersection obeys  convexity condition \ref{con1} above. (The reason we required that $U\subset 2U_{i}$ was so that $g_{i}$ would be defined and convex on all of $U$.) 

The condition that the transition map between the coordinate chart $U$ and the coordinate chart where $g_{i}$ satisfies condition \ref{con2} is equivariant, and the fact that the tropical part of $U$ is some face of this coordinate chart implies that we can use Remark \ref{pullback monomial} to see that $g_{i}$ in the coordinate chart $\mathbb R^{n}\times\et kP$ containing $U$ is some sum of positive functions of $\mathbb R^{n}$ times the square absolute value of smooth monomial functions. Therefore convexity condition \ref{con2} also holds for any intersection of our $U$'s.

We may continue covering strata of lower tropical dimension until a finite cover of $\ex B$ satisfying Claim \ref{con cover} has been constructed.

\begin{claim}\label{con good}If  $U\subset \mathbb R^{n}\times \et kP$ satisfies the convexity conditions \ref{con1} and \ref{con2} of Claim \ref{con cover}, then $U$ is  isomorphic to $\mathbb R^{n}\times \et kP$.
\end{claim}

To prove Claim \ref{con good}, note that the set $U'$ defined by $U\cap \mathbb R^{n}\times \et k{P^{\circ}}=U'\times \et k{P^{\circ}}$ is geodesically convex and open, and therefore diffeomorphic to $\mathbb R^{n}$. We can therefore reduce to the case that $U'=\mathbb R^{n}$.   

Recall that $U$ is equal to the set where $g_{i}<1$ for some finite number of functions $g_{i}$ which are sums of positive functions on $\mathbb R^{n}$ times square absolute values of monomial functions. Suppose that $g_{i}=g+y$ where $g$ and $y$ are nonnegative and $y$ is (the pull back of) a function on $\mathbb R^{n}$. We may replace such a $g_{i}$ by  $g(1-y)^{-1}$. Therefore, we may assume that these $g_{i}$ are a sum of positive functions times the square absolute value of nonconstant monomial functions.  Choose a diffeomorphism $\rho\co [0,1)\longrightarrow [0,\infty)$ so that close to $0$, $\rho(x)=x$. The function 
\[G:=\sum_{i}\rho\circ g_{i}\] is smooth and proper on each $\et mP$ fiber of $U$.  In what follows, we shall use $G$ to define a complete vector field $v$ on $U$ whose (negative time) flow eventually sends each point of $U$ into an arbitrarily small neighborhood of $U'\times \et k{P^{\circ}}$ where $G$ is small. We shall also define an analogous complete vectorfield $v'$ on $\mathbb R^{n}\times \et mP$, so that $v'$ agrees with $v$ when $G$ is small. Then we shall define an isomorphism $U\longrightarrow \mathbb R^{n}\times \et mP$ by a map locally defined by flowing $v$ for some time $-T$ so that $G$ becomes small, then flowing back out for time  $T$ using $v'$.

Note that at each point of $\et mP$, there exists a vector $v$ so that for all smooth monomials $\zeta$,  $v\abs\zeta^{2}$ is positive if  $\zeta$ is nonzero. Therefore, $vG>0$ if $G\neq 0$. Therefore, if $\nabla G$ indicates the gradient of $G$ in the $\et mP$ direction using the standard flat metric, then $\nabla G$ is nonzero whenever $G$ is nonzero. Let $v$ be a smooth vectorfield on $U$ so that \begin{itemize}\item$v$ is in the kernel of the projection  $U\longrightarrow \mathbb R^{n}$, \item $0\leq vG \leq 1$ \item and $vG>0$ when $G>0$. \end{itemize} This vectorfield $v$ is complete on $U$ and for any point $p\in U$ and $\epsilon>0$, there exists some time $T$ so that $G(\Phi_{-tv}(p))<\epsilon$ for all $t>T$. 

We shall now define an analogous vectorfield $v'$ on $\mathbb R^{n}\times \et mP$. Let
\[G':=\sum_{i}g_{i}\]
 When $G$ is small,  $G= G'=\sum_{i}g_{i}$. Note that $G'$ is proper restricted to $\et mP$ fibers of $\mathbb R^{n}\times\et mP$ and $\abs{\nabla G'}>0$ whenever $G'>0$. We can therefore choose some smooth vectorfield $v'$ so that\begin{itemize}\item $v'=v$  on a neighborhood of $\et m{P^{\circ}}$,\item $v'$ is in the kernel of the projection to $\mathbb R^{n}$,\item $v'G'>0$ wherever $G'>0$,\item and $v' G'\leq 1$. \end{itemize}

Consider the map $U\longrightarrow \mathbb R^{n}\times \et mP$ given by the limit as $t\to\infty$ of   \[\Phi_{tv'}\circ\Phi_{-tv}\] the flow for time $-t$ of $v$ followed by the flow for time $t$ of $v'$. Note that  $v'$ and $v$ agree in a neighborhood of $\mathbb R^{n}\times\et m{P^{\circ}}$ and $\Phi_{-tv}$ eventually brings any point into this neighborhood. Therefore, around any point, this limit is simply given by   $\Phi_{tv'}\circ\Phi_{-tv}$ for some large $t$. It follows that this map is smooth. It is also obviously invertible, as $\Phi_{-tv'}$ also eventually brings each point into a neighborhood of $\mathbb R^{n}\times \et m{P^{\circ}}$.  It follows that $U$ is isomorphic to $\mathbb R^{n}\times\et mP$, and we have completed the proof Claim \ref{con good}, and therefore the proof of our lemma.

\stop

\section{Cohomology does not change in connected families}

\label{invariant section}

Consider a family of exploded manifolds over $\mathbb R$. Any such family is trivial, so the cohomology of different fibers may be identified. This identification does not depend on the choice of trivialisation, so it is canonical. Similarly, for any family of exploded manifolds over an exploded manifold $\ex F$, given any two points $p_{1}, p_{2}$ in $\ex F$, and a smooth path $\gamma\co [0,1]\longrightarrow \ex F$ from $p_{1}$ to $p_{2}$, there is an identification of the cohomology of the fiber over $p_{1}$ with the fiber over $p_{2}$. This identification only depends on the isotopy class of $\gamma$. Of course, even if $\ex F$ is connected, if $p_{1}$ and $p_{2}$ have different image in $\totb {\ex F}$, there will be no such path. 

On the other hand, consider a family of  exploded manifolds over $\et mP$ where the polytope $P$ is open. A differential form on some fiber $\ex B$ of this family may be regarded as some section of $\totl{\bigwedge T^{*}\ex B}$ over $\totl{\ex B}$. The smooth part of any fiber of this family is canonically isomorphic to $\totl{\ex B}$, and the smooth part of the cotangent space of every fiber is canonically isomorphic to $\totl{T^{*}\ex B}$, therefore there is a canonical identification of differential forms on any fiber of our family with differential forms on $\totl{\ex B}$. This canonical identification preserves exterior differentiation, wedge products and integration, and gives a canonical identification of the cohomology of every fiber of our family. 

Here is a surprising observation: Given a path $\gamma$ between two points in $\et mP$, the identification of cohomology from $\gamma$ will sometimes be different from the canonical identification coming from identifying the smooth parts of every fiber; these two different identifications will always agree if $\gamma$ is a path that only travels in the directions spanned by integral vectors, but these identifications will often be different if the integral of the imaginary part of $\tilde z_{i}^{-1}d\tilde z_{i}$ does not vanish on $\gamma$. 

 We shall take the view that the identification of cohomology coming from a path is the `correct' identification. To identify the cohomology of two fibers with different tropical parts, we shall need the following notion of a long path.

\begin{defn}[long path] Let $\ex X$ be some refinement of $\ex T$. Consider $\ex X$ as a $\Cs\infty1$ exploded manifold. 

Recall from section $3$ of \cite{iec} that  any subset of $\ex X$ inherits the structure of an abstract exploded space. A long line $\ex L$ is an abstract exploded space isomorphic to the subset of $\ex X$ given by 
\[\ex L:=\{\tilde z\in (0,\infty)\e{\mathbb R}\subset \mathbb C^{*}\e{\mathbb R}\}\subset\ex X\]
where $\tilde z$ is the standard coordinate on $\ex T$ pulled back to $\ex X$. The ends of $\ex L$ are the connected components  of $\ex L$ intersected with the strata of $\ex X$ with unbounded tropical part. 

\

A long path in an exploded manifold $\ex B$ is a morphism
\[\gamma\co \ex L\longrightarrow \ex B\]
 of abstract exploded spaces so that $\gamma$ is constant  on a neighborhood of the ends of $\ex L$.

\end{defn}

Expanding this definition a little, $\ex L$ is locally isomorphic either to $\mathbb R$ or a subset of $\et 1{I}$ defined by $\tilde z\in (0,\infty)\e{\mathbb R}$. Our long path $\gamma$  will be an ordinary smooth path in $\ex B$ where $\ex L$ is isomorphic to $\mathbb R$, and on the other coordinate charts it can be given by specifying a $\Cs\infty1$ map $\hat\gamma:\et 1{I}\longrightarrow \ex B$, and restricting to the subset where $\tilde z\in (0,\infty)\e{\mathbb R}$. 
There is a long path joining any two points in a connected  exploded manifold.

A long line $\ex L$ is a type of exploded space which can be regarded as a $\Cs\infty1$ version of an exploded fibration, defined in \cite{gokova}. Given a family $\hat{\ex B}\longrightarrow \ex B_{0}$ of exploded manifolds and a long path $\gamma$ in $\ex B_{0}$, the pullback $\gamma^{*}\hat{\ex B}$ of this family is a nice exploded space on which differential forms may be defined. 
{\setlength\mathsurround{0pt}\[\begin{tikzcd}\gamma^{*}\hat{\ex B}\rar \dar&\hat{\ex B}\dar
\\ \ex L\rar{\gamma}&\ex B_{0}\end{tikzcd}\]}
Where $\ex L$ is isomorphic to $\mathbb R$, $\gamma^{*}\hat{\ex B}$ is isomorphic to a family of exploded manifolds over $\mathbb R$ so differential forms can be defined as usual. Where $\ex L$ is isomorphic to the subset of $\et 1I$ given by $\tilde z\in(0,\infty)\e{\mathbb R}$, $\gamma^{*}\hat{\ex B}$ is an abstract exploded space isomorphic to the subset of some family $\hat\gamma^{*}\hat{\ex B}\longrightarrow \et 1I$ given by restricting to $\tilde z\in(0,\infty)\e{\mathbb R}$.  As in section $6$ of  \cite{iec}, we can define the tangent sheaf of $\gamma^{*}\hat{\ex B}$ as the sheaf of derivations on exploded functions. Where $\ex L$ is isomorphic to $\{\tilde z\in(0,\infty)\e{\mathbb R}\}\subset \et 1I$, the tangent space $T\gamma^{*}\hat{\ex B}$ may be given by restricting the tangent space of $\hat\gamma^{*}\hat{\ex B}$ to  $\gamma^{*}\hat{\ex B}\subset \hat\gamma^{*}\hat{\ex B}$ in the obvious way: here, $T\gamma^{*}\hat{\ex B}$ is  the kernel of the imaginary part of $\tilde z^{-1}dz$  within $T\hat \gamma^{*}\hat{\ex B}\rvert_{\gamma^{*}\hat{\ex B}}$. Integral vectors within $T\gamma^{*}\hat{\ex B}$ are also defined as usual, and correspond to the restriction of integral vectors from $T\hat \gamma^{*}\hat{\ex B}$.

 We can then define $T^{*}\gamma^{*}\hat{\ex B}$ as the dual vectorbundle to $T\gamma^{*}\hat{\ex B}$, and define differential forms on $\gamma^{*}\hat{\ex B}$ as $\Cs\infty1$ sections of $\bigwedge T^{*}\gamma^{*}\hat{\ex B}$. Define $\Omega^{*}(\gamma^{*}\hat{\ex B})$ to be the sheaf of differential forms on $\gamma^{*}\hat{\ex B}$ that are in $\Omega^{*}$ restricted to any fiber of $\gamma^{*}\hat{\ex B}\longrightarrow \ex L$, and that vanish on any integral vectors. Again, $\Omega^{*}(\gamma^{*}\hat{\ex B})$ is as usual where $\ex L$ is isomorphic to $\mathbb R$, and where $\ex L$ is isomorphic to $\{\tilde z\in(0,\infty)\e{\mathbb R}\}\subset \et 1I$, $\Omega^{*}(\gamma^{*}\hat{\ex B})$ consists of differential forms that are locally equal to the restriction of forms in $\Omega^{*}(\hat \gamma^{*}\hat{\ex B})$ to $\gamma^{*}\hat{\ex B}\subset \hat\gamma^{*}\hat{\ex B}$.

\begin{prop}\label{long path extension} Suppose that $\hat{\ex B}\longrightarrow \ex B_{0}$ is a family of exploded manifolds,  $\gamma\co \ex L\longrightarrow \ex B_{0}$ is a long path, and $\ex B$   is a fiber of $\hat{\ex B}\longrightarrow \ex B_{0}$ over one of the ends of $\gamma$. Then any closed differential form  $\theta\in \Omega^{*}\ex B$ extends to a closed differential form $\hat \theta\in\Omega^{*}(\gamma^{*}\hat{\ex B})$. 
\end{prop}

We shall delay the slightly technical proof of Proposition \ref{long path extension} until after we have defined our isomorphism between the cohomology of different fibers of a family of exploded manifolds. 

\begin{defn}\label{Psi gamma} Suppose that $\hat{\ex B}\longrightarrow \ex B_{0}$ is a family of exploded manifolds,  $\gamma\co \ex L\longrightarrow \ex B_{0}$ is a long path, and $\ex B$ and $\ex B'$  are the fibers of $\hat{\ex B}\longrightarrow \ex B_{0}$ over the ends of $\gamma$. Define a map 
\[\Psi_{\gamma}\co H^{*}(\ex B)\longrightarrow H^{*}(\ex B')\]
as follows: Choose a representative $\theta\in\Omega^{*}(\ex B)$  for a given cohomology class $[\theta]\in H^{*}(\ex B)$. Extend $\theta$ to a closed form $\hat \theta\in\Omega^{*}(\gamma^{*}\hat{\ex B})$, then let $\theta'\in\Omega^{*}(\ex B')$ be the restriction of $\hat \theta$ to $\ex B'$. Define $\Psi_{\gamma}([\theta])=[\theta']\in H^{*}(\ex B')$.
\end{defn}

The proposition below tells us that $\Psi_{\gamma}$ is an isomorphism. In order to use integration along the fiber, we have added the assumption that $\hat{\ex B}\longrightarrow\ex B_{0}$ is a family of oriented exploded manifolds. This assumption may be easily removed by using differential forms with coefficients twisted by the orientation bundle of the fibers. To use Poincar\'e duality, we have added the assumption that our exploded manifolds are simplicial at infinity; see Definition \ref{simplicial}.

\begin{prop}\label{invariance prop}If $\hat{\ex B}\longrightarrow \ex B_{0}$ is a family of oriented exploded manifolds that are simplicial at infinity, then the map $\Psi_{\gamma}\co H^{*}(\ex B)\longrightarrow H^{*}(\ex B')$ from Definition \ref{Psi gamma} does not depend on the choice of $\theta$ and $\hat \theta'$. Morover, $\Psi_{\gamma}$ is a linear isomorphism that preserves wedge products and integration. \end{prop}

\pf

We shall begin with the following claim that tells us $\Psi_{\gamma}$ is compatible with integration.
\begin{claim}\label{constant integral} Suppose that $\hat \theta\in \Omega^{*}(\gamma^{*}\hat{\ex B})$ is closed. Then the integral of  $\hat \theta$ over any two fibers of $\gamma^{*}\hat{\ex B}$ is equal.
\end{claim}

To prove Claim \ref{constant integral}, we may assume that $\hat\theta$ has the correct degree to be integrated on the fibers of $\gamma^{*}\hat{\ex B}$. In regions where $\ex L$ is isomorphic to $\mathbb R$, Theorem \ref{fiber integration} implies that  the integral of $\hat{\theta}$ along the fiber of the map $\gamma^{*}\hat{\ex B}\longrightarrow \ex L$ is a constant function, so on such regions, the integral of $\hat{\theta}$ on  fibers of $\gamma^{*}\hat{\ex B}$
is constant. To prove Claim \ref{constant integral}, we need only show that the same holds in regions where $\ex L$ is isomorphic to $\{\tilde z\in(0,\infty)\e{\mathbb R}\}\subset \et 1I$.

On such a region, $\gamma$ is the restriction of some map $\hat\gamma\co \et 1I\longrightarrow \ex B_{0}$, and $\gamma^{*}\hat{\ex B}$ may be regarded as a subset of the family of exploded manifolds $\pi\co \hat\gamma^{*}\hat{\ex B}\longrightarrow \et 1I$. Theorem \ref{fiber integration} tells us that integration along the fiber of $\pi$ is well defined and commutes with exterior differentiation. On this region, there exists a (not necessarily closed) differential form $\hat \theta'$ defined on $\hat\gamma^{*}\hat{\ex B}$ so that $\theta$ is the restriction of $\hat \theta'$ to $\gamma^{*}\hat{\ex B}\subset \hat\gamma^{*}\hat{\ex B}$. Although $\hat \theta'$ need not be closed, the fact that $\hat \theta$ is closed implies that $d\hat\theta'$ restricted to $\gamma^{*}\hat{\ex B}\subset \hat\gamma^{*}\hat{\ex B}$ is $0$.

Theorem \ref{fiber integration} tells us that $\pi_{!}(\hat \theta')$ is a function defined on $\et 1I$, and that $d\pi_{!}(\hat\theta')=\pi_{!}(d\hat\theta')$. As $d\hat\theta'$ vanishes on $\gamma^{*}\hat{\ex B}\subset \hat\gamma^{*}\hat{\ex B}$, it follows that $\pi_{!}(\hat\theta')$ is constant on $\{\tilde z\in(0,\infty)\e{\mathbb R}\}\subset \et 1I$. In other words, the integral of $\hat \theta$ over fibers of $\gamma^{*}\hat{\ex B}$ is constant. This completes the proof of Claim \ref{constant integral}. Note that Claim \ref{constant integral} implies that $\Psi_{\gamma}$ is compatible with integration. 

\

We shall now show that $\Psi_{\gamma}$ is well defined by showing that it sends forms $\theta$ that are $0$ in homology to forms that are zero in homology. In order to use the version of Poincare duality given in Theorem \ref{Poincare duality}, we shall now assume that every map  $\ex T\longrightarrow \ex B$ must be constant. (So long as $\ex B_{0}$ is connected, the same must hold for every fiber of $\hat{\ex B}\longrightarrow \ex B_{0}$.) Suppose that $\hat\theta\in\Omega^{*}(\gamma^{*}\hat{\ex B})$ is closed and restricts to be $0$ in $H^{*}(\ex B)$. Then given any other closed form $\hat\theta'\in \Omega^{*}(\gamma^{*}\hat{\ex B})$, the integral of $\hat \theta\wedge\hat\theta'$ over $\ex B$ must be $0$, therefore Claim \ref{constant integral} implies that the integral of $\hat\theta\wedge\hat\theta'$ over $\ex B'$ must also be $0$. Theorem \ref{Poincare duality} implies that if $\hat\theta$ restricted to $\ex B'$ was not $0$ in $H^{*}(\hat{\ex B})$, then there would exist a closed differential form $\theta'\in\Omega^{*}(\ex B')$ so that the integral of $\hat\theta\wedge \theta'$ over $\ex B'$ is nonzero. Proposition \ref{long path extension} implies that such a $\theta'$ extends to some closed $\hat \theta'\in \Omega^{*}(\gamma^{*}\hat{\ex B})$. Therefore $\hat\theta$ restricts to be $0$ in $H^{*}(\ex B')$ if and only if  it restricts to be $0$ in $H^{*}(\ex B)$.

We have shown under the assumption that all maps $\ex T\longrightarrow \ex B$ are constant and $\ex B$ is simplicial at infinity,  $\Psi_{\gamma}\co H^{*}(\ex B)\longrightarrow H^{*}(\ex B')$ is well defined independent of the choice of $\theta$ and $\hat \theta$ in Definition \ref{Psi gamma}.   In the more general case that there are nonconstant maps $\ex T\longrightarrow \ex B$, then $\hat{\ex B}\longrightarrow \ex B_{0}$ is isomorphic to some $\ex T^{n}$ bundle over some other bundle $\hat {\ex B}'\longrightarrow \ex B_{0}$ whose fibers are still simplicial at infinity but do not admit nonconstant maps from $\ex T$. Any differential form in $\Omega^{*}$ of  a fiber of $\hat{\ex B}\longrightarrow \ex B_{0}$ is equal to the pullback of a unique form on the corresponding fiber of $\hat{\ex B}'\longrightarrow \ex B_{0}$. Similarly, forms in $\gamma^{*}\hat{\ex B}$ are equal to the pullback of forms in $\gamma^{*}\hat{\ex B}'$. Therefore, the fact that $\Psi_{\gamma}\co H^{*}(\ex B)\longrightarrow H^{*}(\ex B')$ is well defined independent of the choice of $\theta$ and $\hat \theta$ follows from the analogous fact for the family $\hat{\ex B}'\longrightarrow \ex B_{0}$.

The fact that $\Psi_{\gamma}$ is well defined independent of choice of $\theta$ and $\hat \theta$ implies that it is linear and preserves wedge products, because $c_{1}\theta_{1}+c_{2}\theta_{2}$ extends to $c_{1}\hat \theta_{1}+c_{2}\hat\theta_{2}$, and $\theta_{1}\wedge\theta_{2}$ extends to $\hat \theta_{1}\wedge \hat\theta_{2}$. Reversing the roles of $\ex B$ and $\ex B'$ gives an inverse to $\Psi_{\gamma}$, so it is a linear isomorphism that preserves wedge products and integration, as required. 

\stop

\

The remainder of this paper is devoted to the  proof of Proposition \ref{long path extension}.

\

 We must extend a closed differential form $\theta\in\Omega^{*}\ex B$ to a closed differential form $\hat\theta\in\Omega^{*}(\gamma^{*}\hat{\ex B})$, where $\ex B$ is the fiber of $\hat{\ex B}\longrightarrow \ex B_{0}$ over an end of $\gamma\co \ex L\longrightarrow \ex B_{0}$.   Because $\gamma$ must be constant on a neighborhood of the ends of $\ex L$, we may extend  $\theta$ to some closed differential form $\hat\theta\in\Omega^{*}(\gamma^{*}\hat{\ex B})$ defined on a neighborhood of $\ex B$.
 
 We can continue to extend $\hat \theta$ to the rest of $\gamma^{*}\hat{\ex B}$ using the following claim that will take several pages to prove:
 \begin{claim}\label{extension} Suppose that an open subset $U$ of $\ex L$ is isomorphic to $\{\tilde z\in(0,\infty)\e{\mathbb R}\}\subset \et 1I$, where $I\subset [0,\infty)$ is some interval containing $0$. Let $\hat\theta$ be a closed differential form in $\Omega^{*}(\gamma^{*}\hat{\ex B})$ defined  where $\abs{\tilde z}>1$. Then there exists a closed differential form in $\Omega^{*}(\gamma^{*}\hat{\ex B})$  defined on all of $U$ that agrees with $\hat\theta$ on the subset where $\abs {\tilde z}>2$. 
 \end{claim}
 
 At this point, we shall use the notion of equivariant coordinate charts, defined and constructed in Appendix A of \cite{cem}. The idea shall be to represent a given homology class in a fiber of $\gamma^{*}\hat{\ex B}$ with a suitably equivariant differential form that extends uniquely as an equivariant differential form on $\gamma^{*}\hat{\ex B}$. 
 
 On our open subset, $\gamma$ is the restriction of a $\Cs\infty1$ map $\hat \gamma\co \et 1I\longrightarrow \ex B_{0}$.
  Lemma A3 of \cite{cem} constructs equivariant coordinate charts for any smooth family of exploded manifolds, however the same proof works to construct equivariant coordinate charts for any $\Cs\infty1$ family of exploded manifolds. The only change required in the proof is that $\Cs\infty1$ vectorfields must be used in place of smooth vectorfields. As in the proof of Lemma A3, we may choose equivariant coordinate charts on $\hat\gamma^{*}\hat{\ex B}$ so that the projection to $\et 1I$ is equivariant. (An alternate proof of the existence of equivariant charts for $\Cs\infty1$ exploded manifolds is contained in the construction of normally rigid structures given in \cite{icc}.)

Given equivariant coordinate charts on $\hat \gamma^{*}\hat{\ex B}$, after shrinking our coordinate charts appropriately,  we may choose  real functions $d_{S}$ on $\hat{\ex B}$ (designed to measure the distance to strata $S$) satisfying the following:
\begin{enumerate}
\item $d_{S}$ is a nonnegative real function whose vanishing set is the closure of $S$.
\item On any of our coordinate charts intersecting $S$ but not the boundary of $S$,  if our chart is (equivariantly) isomorphic to an open subset of $\mathbb R^{n}\times \et mP$, then
\[d_{S}=\sum_{i}f_{i} \abs{\zeta_{i}}^{2}\]
where
\begin{enumerate}\item each $\zeta_{i}$ is a smooth monomial on $\et mP$.
\item $f_{i}$ is a nonnegative smooth function on $\mathbb R^{n}$.
\end{enumerate}
\item \label{cond3} If such a coordinate chart intersects a stratum $T$, then on this coordinate chart, 
\[d_{T}=\Delta_{T}d_{S}\]
 In other words, if $d_{S}=\sum _{i}f_{i}\abs{\zeta_{i}}^{2}$, then $d_{T}$ is the corresponding sum in which every monomial $\zeta_{i}$ that does not vanish on $T$ has been removed. 
\end{enumerate}

 Such distance functions may be constructed starting with the strata with the largest dimensional tropical part, then choosing the distance functions for the smaller strata compatibly after shrinking coordinate charts a little. Once such distance functions are chosen, there exist distances $R_{S}>0$ and an open covering of  $\hat \gamma^{*}\hat{\ex B}$ by a new collection of equivariant coordinate charts $\hat U$ compatible with our old equivariant coordinate charts, and satisfying the following conditions.
 \begin{enumerate}
 \item\label{u1}  Each $\hat U$ that intersects the stratum $S$ but not the boundary of $S$ is equivariantly isomorphic to an open subset of $\mathbb R^{n}\times \et mP$ where $d_{S}<2R_{S}$.
  \item\label{u2} The subset where $d_{S}<2R_{S}$ is covered by charts that intersect $S$.
 \item\label{u3} For any stratum $T$, the subset where $d_{T}>0$ is covered by charts that do not intersect $T$.
 \item \label{u4}If $T\neq S$ is in the closure of $S$, then $R_{T}>R_{S}$.
 \end{enumerate}

Choose some fiber $\ex B$ of $\gamma^{*}\hat{\ex B}$ over the subset of $\et 1I$ where $1>\abs{\tilde z}>0$ and close enough to the central stratum of $\et 1I$ to be covered by the subsets where $d_{S}<R_{S}/2$ for strata $S$ that project to the central stratum of $\et 1I$.

In the pages that follow, we shall show that we can replace a given closed differential form on $\ex B$ with a differential form that is equivariant where $d_{S}<R_{S}/2$ in the following sense:  Suppose that $\hat U$ is one of our charts intersecting $S$. Then $\hat U$ is equivariantly isomorphic to an open subset of $\mathbb R^{n}\times \et mP$. Suppose that the stratum $S$ correponds to a $k$--dimensional stratum of $P$ that is the interior of the face of $P$ defined by setting the last $m-k$ coordinates on $P\subset \mathbb R^{m}$ equal to $0$. Consider the lattice $\mathcal I_{S}$  of vector fields consisting of integer sums of  the real parts of $\{\tilde z_{1}\frac \partial{\partial z_{1}},\dotsc,\tilde z_{k}\frac \partial{\partial z_{k}}\}$, and the corresponding lattice $\mathcal I_{S}'$ consisting of integer sums of the imaginary parts of $\{\tilde z_{1}\frac \partial{\partial z_{1}},\dotsc,\tilde z_{k}\frac \partial{\partial z_{k}}\}$.  Any equivariant change of coordinate charts intersecting $S$ preserves these lattices $\mathcal I_{S}$ and $\mathcal I'_{S}$, so $\mathcal I_{S}$ and $\mathcal I'_{S}$ are well defined sheaves of vectorfields on $\hat\gamma^{*}\hat{\ex B}$ where $d_{S}<2R_{S}$. The restriction $U$ of $\hat U$ to $\ex B$ is an open subset of $\mathbb R^{n'}\times (\mathbb C^{*})^{l}\times \et {m'}{P'}$, and the restriction of the vector fields in $\mathcal I_{S}$ and $\mathcal I_{S}'$ now include the  real and imaginary part of $z_{i}\frac \partial{\partial z_{i}}$ where $z_{i}$ indicates the coordinate on the $i$th $\mathbb C^{*}$. 

Say that a differential form $\theta$  is $S$--equivariant if $L_{v}\theta=0$ for all $v$ in $\mathcal I_{S}\cup\mathcal I_{S}'$ and $i_{v}\theta=0$ for all $v\in \mathcal I_{S}$. We shall also use the weaker condition that $\theta$ is $\mathcal I'_{S}$--equivariant if $L_{v}\theta=0$ for all $v\in \mathcal I_{S}'$. Say that a function or a subset of a coordinate chart is $\mathcal I'_{S}$--invariant if it is preserved by the flow of vectorfields in $\mathcal I'_{S}$. In what follows, we wish to replace a given closed differential form with an equivariant differential form representing the same homology class. We shall do this by first replacing it with a  $\mathcal I_{S}'$--equivariant form, then modifying this form to achieve $S$--equivariance.
 
 \begin{claim}\label{T equivariant} Every cohomology class in  $H^{*}\ex B$ is represented by some closed form $\theta\in \Omega^{*}(\ex B)$ that is $\mathcal I'_{S}$--equivariant where $d_{S}<R_{S}$ for all strata $S$ of $\hat\gamma^{*}\hat{\ex B}$.
 
 \end{claim}
 
Let $U$ be some $\mathcal I'_{S}$--invariant subset of $\ex B$ where $\mathcal I_{S}'$ makes sense, and on which the sheaf $\mathcal I'_{S}$ is generated by global sections.
  Choose a basis $\{v_{l}\}$ for $\mathcal I'_{S}$ on $U$, then define
\[K_{l}(\theta):=\frac 1{2\pi}\int_{0}^{1}\int_{0}^{2\pi}-\Phi^{*}_{stv_{l}}i_{tv_{l}}\theta ds dt\] 
where  $\Phi_{st v_{l}}$ indicates the flow of $v_{l}$ for time $st$. This flow is complete on $U$ or any $\mathcal I'_{S}$--invariant subset, so $K_{l}(\theta)$ is well defined for $\theta\in \Omega^{*}(U)$, and also makes sense on any $\mathcal I'_{S}$--invariant subset. Note that if on some $\mathcal I_{S}'$--invariant subset,  $L_{v}\theta=0$ for any $v\in\mathcal I'_{S}$,   then $L_{v}K_{l}\theta=0$ on this subset too. The following calculation shows how $K_{l}$ can be used to modify a given form to a $v_{l}$--invariant form:
\[\begin{split}(dK_{l}+K_{l}d)\theta &=\frac 1{2\pi}\int_{0}^{1}\int_{0}^{2\pi}-\Phi^{*}_{stv_{l}}L_{tv_{l}}\theta ds dt
\\ &= \frac 1{2\pi}\int_{0}^{2\pi}\theta-\Phi^{*}_{tv_{l}}\theta dt
\\ &=\theta-\frac 1{2\pi}\int_{0}^{2\pi}\Phi^{*}_{tv_{l}}\theta dt\end{split}\]
Choose a $\mathcal I'_{S}$--invariant subset $V$ of $U$, and an $\mathcal I'_{S}$--invariant cut-off function $\rho$ that is $1$ restricted to $V$, and has compact support within $U$. Then for any closed form $\theta\in \Omega^{*}(U)$, we may replace $\theta$ with
\[\prod_{l}(1-d\rho K_{l})\theta\]
This modified form represents the same cohomology class, but is now $\mathcal I_{S}'$--equivariant on $V$. Moreover, if $\mathcal I'_{T}\subseteq \mathcal I'_{S}$, this new form remains $\mathcal I'_{T}$--equivariant on any $\mathcal I'_{S}$--invariant subset on which $\theta$ was already $\mathcal I'_{T}$--equivariant. Condition (\ref{cond3}) on $d_{T}$ implies that it is $\mathcal I'_{S}$--invariant wherever $\mathcal I'_{S}$ makes sense. If $\theta$ is already $\mathcal I'_{T}$--equivariant where $d_{T}<R_{T}$, for any collection of strata $T$ not in the closure of $S$,  we may cover the region where $d_{S}<R_{S}$ with $\mathcal I'_{S}$--equivariant subsets $V$ as above, and modify $\theta$ using the above procedure to a form that is $\mathcal I'_{S}$--equivariant where $d_{S}<R_{S}$, and still $\mathcal I'_{T}$--equivariant where $d_{T}<R_{T}$.   Claim \ref{T equivariant} then follows by induction on the tropical dimension of the strata involved.

\

Consider the charts $\hat U$ satisfying the conditions enumerated on page \pageref{u1}. For each such  chart $\hat U$ intersecting $S$ but not any other strata in the closure of $S$, choose  $\hat U'$ to be a compactly contained open subset of $\hat U$ where each $T^{m}_{P}$ fiber is either empty or the subset where $d_{S}<R_{S}/2$.  Choose the charts $\hat U'$ to still cover $\hat\gamma^{*}\hat{\ex B}$ so that the subset where $d_{S}<R_{S}/2$ is covered by charts  $\hat U'$ that intersect $S$. (The conditions enumerated for the charts $\hat U$ ensure that this is possible.)  Let $U$ be the intersection of $\hat U$ with $\ex B$ and $U'$ be the intersection of $\hat U'$ with $\ex B$.

\

 \begin{claim}\label{equivariant modification}Let $U_{i}$ be one of our coordinate charts and let $S$ be the stratum intersecting $\hat U_{i}$ with maximal tropical dimension. Let $\theta$ be any closed differential form  in $\Omega^{*}(U_{i})$ so that
 \begin{itemize}\item $\theta$ is  $\mathcal I_{T}'$--equivariant  wherever $d_{T}<R_{T}$ and $S$ is in the closure of $T$,
 \item  and $\theta$ is also  $T'$--equivariant where $d_{T'}<R_{T'}/2$ whenever $T'\neq S$ is in the closure of $S$. 
  \end{itemize}
   Then there exists a closed differential form $\theta' \in \Omega^{*}(U_{i})$ representing the same cohomology class  so that: \begin{itemize}\item $\theta'$ is $S$--equivariant on $U'_{i}\subset U_{i}$. 
 \item The support of $\theta-\theta'$ within $U_{i}$ is compact and contained in the set where $d_{S}<R_{S}$.
 \item If $\theta$ is already $S$--equivariant on $U_{i}\cap U'_{j}$ and $\hat U'_{j}$ intersects $S$ then  $\theta'=\theta$  on $U_{i}\cap U'_{j}$.
 \item  Where $d_{T}<R_{T}$, and $S$ is in the closure of $T$, $\theta'$ is $\mathcal I_{T}'$--equivariant.
 \item  Where  $d_{T}<R_{T}/2$ and $T\neq S$ is in the closure of $S$, $\theta'$ is $T$--equivariant. \end{itemize}
 \end{claim}

 We have that $U_{i}$ is an open subset of $\mathbb R^{n}\times (\mathbb C^{*})^{k}\times \et mP$. Define $K_{0}$ to be the operator that restricted to each $\et mP$ fiber is the $K$ from equation (\ref{defK}) on page \pageref{defK}. Because smooth monomials increase along flowlines of the vectorfield $v$ from equation (\ref{v def}) used to define $K$, the negative time flow of $v$ decreases $d_{S}$, and condition (\ref{cond3}) on $d_{T}$ implies that the negative time flow of $v$ also decreases $d_{T}$ where $d_{T}<2R_{T}$ for all other strata $T$. Therefore such a $K_{0}$ is well defined on $U_{i}$, the subset of $U_{i}$ where $d_{T}<R_{T}$ or $R_{T}/2$,  and also on $U_{i}\cap U_{j}'$ so long as $\hat U_{j}$ intersects $S$.  To see that $K_{0}\theta$ is in $\Omega^{*}(U_{i})$, note that the estimate (\ref{Kestimate}) on page \pageref{Kestimate} now  holds with the constant $c$ replaced by a continuous function, and the rest of the proof that $K(\theta)$ is in $\Omega^{*}$ applies to show that $K_{0}(\theta)$ is in $\Omega^{*}(U_{i})$. 
   Because $v$ commutes with the vectorfields in $\mathcal I'_{S}$, $K_{0}(\theta)$, like $\theta$, is $\mathcal I'_{S}$--equivariant where $d_{S}<R_{S}$.  Because  $v$ is in the span of the vectorfields  $\mathcal I_{S}$,  $K_{0}(\theta)$ vanishes on the set $U_{i}\cap U'_{j}$ whenever $\theta$ is $S$--equivariant on $U_{i}\cap U'_{j}$ (and $\hat U_{j}$ intersects $S$). 
Similarly $K_{0}(\theta)$ vanishes where $d_{T}<R_{T}/2$  and $T\neq S$ is in the closure of $S$. Equation (\ref{Kprop}) on page \pageref{Kprop} also applies to $K_{0}$, giving that $\theta-(K_{0}d+dK_{0})\theta$ is independent of the $\et mP$ coordinates. 

We also need to modify $\theta$ to vanish in the directions given by multiplying the $(\mathbb C^{*})^{k}$ coordinates by real numbers. Let $x_{l}:=\log\abs{z_{l}}$, where $z_{l}$ is the $l$th coordinate on $(\mathbb C^{*})^{k}$. Our distance function $d_{S}$  restricted to a $(\mathbb C^{*})^{k}$ fiber is in the form $\sum c_{j}e^{2\alpha^{j}\cdot x}$  where  $c_{j}$ are positive real numbers and  $\alpha^{j}$ are integral vectors whose positive span contains $\mathbb R^{k}$. In particular, on such a fiber, $d_{S}$ is a strictly convex function of $x$, whose double derivative (in $x$) is always positive definite. In particular, this implies that $d_{S}$ achieves a unique minimum on each $(\mathbb C^{*})^{k}$ fiber, and that the position where this minimum is achieved depends depends smoothly on the other $\mathbb R^{n}\times \et mP$ coordinates. Let 
\[H_{t}\co [0,1]\times\mathbb R^{n}\times(\mathbb C^{*})^{k}\times \et mP\longrightarrow \mathbb R^{n}\times(\mathbb C^{*})^{k}\times \et mP\]
be the homotopy that is the identity on the $\mathbb R^{n}$, $\et mP$, and the angular $(\mathbb C^{*})^{k}$ coordinates, and that is linear in the $x$ coordinates, and that at $t=0$ is the identity, and at $t=1$ sends the  $x$ coordinate on each fiber to the unique value at which $d_{S}$ is minimised.  Then define
\[K'\theta:=\int_{0}^{1}-H_{t}^{*}i_{\frac {\partial H_{t}}{\partial t}}\theta dt\]
\[(dK'+K'd)\theta=\int_{0}^{1}-H_{t}^{*}L_{\frac {\partial H_{t}}{\partial t}}\theta dt=\theta-H_{1}^{*}\theta\]
Note that $\frac{\partial H_{t}}{\partial t}$ is in the  span of the vectorsfields in $\mathcal I_{S}$. If $\hat U_{j}$ also intersects $S$,  the intersection of $U'_{j}$ with a $(\mathbb C^{*})^{k}$ fiber is either empty, or the subset where $d_{T}<R_{T}/2$ for some $T$ in the closure of $S$. Condition (\ref{cond3}) on $d_{T}$ implies that restricted to such a fiber, $d_{T}$ is equal to $d_{S}$ plus a constant  therefore our homotopy sends $U_{i}\cap U_{j}'$ into itself.  It follows that if $\hat U_{j}$ intersects $S$, and $\theta$ is already $S$--equivariant on $U_{i}\cap U_{j}'$, then  $K'(\theta)$ vanishes on $U_{i}\cap U_{j}'$. Similarly, $K'(\theta)$ vanishes on the set where $d_{T}<R_{T}/2$ for any $T\neq S$ in the closure of $S$. As $H_{t}$ sends the set where $d_{S}<R_{S}$ inside itself,  the observation that $\frac{\partial H_{t}}{\partial t}$ commutes with the vectorfields  in $\mathcal I'_{S}$ implies that $K'(\theta)$ is $\mathcal I'_{S}$--equivariant where $d_{S}<R_{S}$.

 Let $\rho$ be some compactly supported, $\mathcal I'_{S}$--invariant smooth function on $U_{i}$ that is $1$ on $U_{i}'$ and $0$ where $d_{S}\geq R_{S}$.  Then for closed $\theta$, define
\[\theta':= (1-d\rho K_{0})(1-d\rho K')\theta\]  
 This $\theta'$ satisfies all the conditions required by Claim \ref{equivariant modification}. In particular, it is $S$--equivariant where $\rho=1$, coincides with $\theta$ where $\rho=0$, also coincides with $\theta$ where $d_{T}<R_{T}/2$ and $T\neq S$ is in the closure of $S$, and on $U_{i}\cap U_{j}'$ if $\theta$ was already $S$--equivariant on $U'_{j}$ and $\hat U_{j}$ intersects $S$.   $\theta'$ is also $\mathcal I'_{S}$--equivariant where $d_{S}<R_{S}$. If  the closure of $T$ contains $S$, being $\mathcal I'_{S}$--equivariant is stronger than being $\mathcal I'_{T}$--equivariant; where $d_{S}<R_{S}$, $\theta'$ is $\mathcal I'_{S}$--equivariant, and therefore $\mathcal I'_{T}$--equivariant, and where $d_{S}\geq R_{S}$, $\theta'=\theta$ so $\theta'$ is therefore $\mathcal I'_{T}$--equivariant where $\theta$ is.

 \
 
 Using Claim \ref{T equivariant}, we may represent any homology class in $H^{*}\ex B$ by a differential form that is $\mathcal I'_{T}$--equivariant where $d_{T}<R_{T}$.  Then starting with our coordinate charts that have maximal tropical dimension, we may use Claim  \ref{equivariant modification} to modify our differential form to a differential form $\theta$ that represents the same cohomology class, but is $T$--equivariant where $d_{T}<R_{T}/2$. In particular, if $S$ has the maximal dimension of the strata that intersect $\hat U$, then $\theta$ is $S$--equivariant on $U'$. Such a differential form has a unique $S$--equivariant extension to the intersection of $\hat U'$ with $ \gamma^{*}\hat{\ex B}\subset \hat\gamma^{*}\hat{\ex B}$.
 
  In particular,  we may sit $\hat U'$ equivariantly inside $\mathbb R^{n}\times \et mP$ so that the first coordinate $\tilde z_{1}$ on $\et mP$ is the projection to $\et 1I$. Let $\theta_{i}$ indicate the imaginary part of $\tilde z_{i}^{-1}d\tilde z_{i}$. Recall that there is some open subset $O\subset \mathbb R^{n}$ so that $\hat U'$ is the subset of $O\times \et mP$ where $d_{S}<R_{S}/2$. The $S$--equivariant differential forms on $\hat U'$ can be written uniquely as $\sum_{w\subset\{1,\dotsc,m\}}\alpha_{w}\bigwedge_{i\in w}\theta_{i}$ where $\alpha_{w}$ is a form pulled back from $O\subset \mathbb R^{n}$. $U'$ corresponds to restricting $\hat U'$ to where $\tilde z_{1}$ is equal to some small real constant. The subset of $\hat U'$ in $\gamma^{*}\hat{\ex B}$ is the subset where $\tilde z_{1}\in(0,\infty)\e{\mathbb R}$. Therefore the kernel of the restriction of such $S$--equivariant forms to $U'$ or the intersection of $\hat U'$ with $\gamma^{*}\hat{\ex B}$ consists of those forms containing $\theta_{1}$. Accordingly, the $S$--equivariant forms on both $ U'$ and the intersection of $\hat U'$ with $\gamma^{*}\hat{\ex B}$ may be written uniquely as $\sum_{w\subset\{2,\dotsc,m\}}\alpha_{w}\bigwedge_{i\in w}\theta_{i}$ where $\alpha_{w}$ is a form pulled back from $O$. 
 
 Let $N$ indicate the subset of $\gamma^{*}\hat{\ex B} $ covered by the charts $\hat U'$ whose projection contains the central strata of $\et 1I$. We have shown that any closed form $\theta\in \Omega^{*}\ex B$ that is $S$--equivariant on $\ex B$ where $d_{S}<R_{S}/2$ extends uniquely to a closed form $\theta'$ in $\Omega^{*}N$ that is $S$--equivariant where $d_{S}<R_{S}/2$. In particular, there exists a closed form $\theta'\in \Omega^{*}N$ that restricts to give a closed form in any chosen cohomology class in $H^{*}\ex B$.  Our neighborhood $N$ contains all the fibers of the family $\gamma^{*}\hat{\ex B}$ in some neighbourhood of the central strata of $\et 1I$. As we may squish all of $\et 1I$ into this neighborhood, it follows that there exists a closed differential form defined on all of $\gamma^{*}\hat{\ex B}$ over $\et 1I$ that restricts to give any chosen cohomology class on the fiber over $\tilde z=2$.    Claim \ref{extension} follows: Given a closed differential form $\hat \theta$ on $\gamma^{*}\hat{\ex B}$ defined over $\et 1I$ where $\tilde z>1$,   we can choose a closed differential form $\hat \theta'$ defined on the subset of $\gamma^{*}\hat{\ex B}$ over $\et 1I$ so that the restriction of $\hat \theta'$ to the fiber where $\tilde z=2$ to gives the same cohomology class as the restriction of $\hat \theta$. As $H^{*}(\ex B)=H^{*}(\ex B\times\mathbb R)$, there exists a differential form $\alpha$ defined on $\gamma^{*}\hat{\ex B}$ so that where $\tilde z<1$, $\hat \theta-\hat\theta'=d\alpha$.  Choose a cutoff function $\rho\co \et 1I\longrightarrow [0,1]$ that is $0$ when $\abs{\tilde z}<1$ and $1$ when $\abs{\tilde z}>2$. Then $\hat \theta'-d\rho\alpha$ is the required extension of $\hat \theta$.  
 
 Claim \ref{extension} allows us to finish the proof of Proposition \ref{long path extension}. As our family $\gamma^{*}\hat{\ex B}$ may be trivialized in a neighborhood of an end of $\ex L$, we may extend a closed differential form defined on the fiber of $\gamma^{*}\hat{\ex B}$ at an end of $\ex L$ to a neighborhood of that end of $\ex L$. The rest of $\gamma^{*}\hat{\ex B}$ is covered by a finite number of coordinate charts on which Claim \ref{extension} allows us to extend our closed differential form. Therefore, our given closed differential form extends to a closed differential form on all of $\gamma^{*}\hat{\ex B}$, as required by Proposition \ref{long path extension}.
 
 \section {Appendix: Exploded manifolds that are symplicial at infinity}
 
 This appendix was added after publication to fix an error in Lemma \ref{local poincare}, which was originally stated without the assumption that $P$ is simplicial at infinity.  The calculation of compactly supported cohomology in Lemma \ref{local poincare} is incorrect without this assumption. An example where it fails is $\et 3P$ where $P$ is a cone over a square. 
 
 As a consequence, the perfect pairing of Theorem \ref{Poincare duality}, and the proof of invariance of cohomology in families of exploded manifolds, Proposition \ref{invariance prop}, both require the extra assumption that our exploded manifold is simplicial at infinity.
 
 \begin{defn}\label{simplicial} A complete polytope $P$ is {\bf standard at infinity} if the monoid of integral vectorfields on $P$ whose flow exists for all positive time is isomorphic to $\mathbb Z^a\times \mathbb N^b$ for some $a$ and $b$.  
 
 More generally, $P$ is {\bf simplicial at infinity} if the monoid of constant vectorfields on $P$ whose flow exists for all positive time is isomorphic to $\mathbb R^a\times[0,\infty)^b$.
 
 An exploded manifold $\ex B$ is {\bf standard or simplicial at infinity} if all polytopes in the tropical part of $\ex B$ are standard or simplicial at infinity.

  \end{defn}
 
As a trivial case, $P$ is standard at infinity if it is compact. 

 Note that the property of being standard or simplicial at infinity is preserved in connected families of exploded manifolds. Roughly speaking, a polytope is standard at infinity if it is a deformation of $\mathbb R^a\times[0,\infty)^b$ and is symplicial at infinity if it is a deformation of $\mathbb R^a$ times a cone that is the non-negative span of a set of $\mathbb R$--linearly independent integral vectors.

\bibliographystyle{plain}
\bibliography{ref.bib}
\end{document}